\documentclass[12pt]{article}

\newbox\squ  
\setbox\squ=\hbox{\vrule width.3pt
            \vbox{\hrule height.3pt width.4em\kern1ex\hrule height.3pt}%
            \vrule width.3pt}
\def\endproof{%
  \ifmmode\eqno\copy\squ\smallskip
  \else{\unskip\nobreak\hfil%
    \penalty50\hskip2em\hbox{}\nobreak\hfil\copy\squ
    \parfillskip=0pt \finalhyphendemerits=0\penalty-100\smallskip}
  \fi}

\usepackage{amsmath}
\usepackage{amssymb}

\begin{document}


\newcommand{\epsi}{\ve}
\newcommand{\M}{{\mathcal M}}
\newcommand{\LL}{{\mathcal L}}
\newcommand{\K}{{\mathcal K}}
\newcommand{\T}{{\mathcal T}}
\newcommand{\A}{{\mathcal A}}
\newcommand{\PP}{{\mathbb P}}
\newcommand{\R}{{\mathbb R}}

\newcommand{\supp}{{\rm{supp}\ts}}
\newcommand{\Card}{{\rm{Card}\ts}}
\newcommand{\dom}{{\rm{dom}\ts}}
\newcommand{\range}{{\rm{range}\ts}}
\newcommand{\rank}{{\rm{rank}\ts}}
\newcommand{\ord}{{\rm{ord}\ts}}
\newcommand{\gr}{{\rm{gr}\ts}}
\newcommand{\const}{{\rm{const}\ts}}
\newcommand{\End}{{\rm{End}\ts}}


\newcommand{\non}{\nonumber}
\newcommand{\wt}{\widetilde}
\newcommand{\wh}{\widehat}
\newcommand{\ot}{\otimes}
\newcommand{\hra}{\hookrightarrow}
\newcommand{\g}{\mathfrak g}
\newcommand{\ve}{\varepsilon}
\newcommand{\ts}{\,}
\newcommand{\U}{ {\rm U}}
\newcommand{\J}{ {\rm J}}
\newcommand{\Y}{ {\rm Y}}
\newcommand{\pr}{{\rm pr}}
\newcommand{\Fun}{{\rm Fun}\ts}
\newcommand{\HH}{{\mathcal H}}
\newcommand{\wH}{\wt{\HH}}
\newcommand{\C}{\mathbb{C}}
\newcommand{\N}{\mathbb{N}}
\newcommand{\Z}{\mathbb{Z}}
\newcommand{\h}{\mathfrak h}
\newcommand{\gl}{\mathfrak{gl}}
\newcommand{\oa}{\mathfrak{o}}
\newcommand{\spa}{\mathfrak{sp}}
\newcommand{\Proof}{\noindent{\it Proof.}\ \ }   
\renewcommand{\theequation}{\arabic{section}.\arabic{equation}} 

\newtheorem{thm}{Theorem}[section]
\newtheorem{prop}[thm]{Proposition}
\newtheorem{cor}[thm]{Corollary}
\newtheorem{defin}[thm]{Definition}
\newtheorem{example}[thm]{Example}
\newtheorem{lem}[thm]{Lemma}

\newcommand{\bth}{\begin{thm}}
\renewcommand{\eth}{\end{thm}}
\newcommand{\bpr}{\begin{prop}}
\newcommand{\epr}{\end{prop}}
\newcommand{\ble}{\begin{lem}}
\newcommand{\ele}{\end{lem}}
\newcommand{\bco}{\begin{cor}}
\newcommand{\eco}{\end{cor}}
\newcommand{\bde}{\begin{defin}}
\newcommand{\ede}{\end{defin}}
\newcommand{\bex}{\begin{example}}
\newcommand{\eex}{\end{example}}

\newcommand{\beq}{\begin{equation}}
\newcommand{\eeq}{\end{equation}}
\newcommand{\bsp}{\begin{split}}
\newcommand{\esp}{\end{split}}

\title{\Large\bf  Degenerate affine Hecke algebras
and centralizer construction for the symmetric groups}
\author{{\sc A. I. Molev\quad and\quad G. I. Olshanski}}


\date{} 
\maketitle

\vspace{7 mm}

\begin{abstract}
In our recent papers
the centralizer construction
was applied to the series of classical Lie algebras
to produce the quantum algebras called
(twisted) Yangians. Here we extend this construction
to the series of the symmetric groups $S(n)$.
We study the `stable' properties of the centralizers
of $S(n-m)$ in the group algebra $\C[S(n)]$
as $n\to\infty$ with $m$ fixed.
We construct a limit centralizer algebra $A$
and describe its algebraic structure. The algebra $A$
turns out to be closely related with
the degenerate affine Hecke algebras.
We also show that the so-called tame representations
of $S(\infty)$ yield a class of natural $A$-modules.

\end{abstract}

\vspace{24 mm}



\noindent
School of Mathematics and Statistics\newline
University of Sydney,
NSW 2006, Australia\newline
alexm@maths.usyd.edu.au

\vspace{7 mm}

\noindent
Dobrushin Mathematics Laboratory\newline
Institute for Problems of Information Transmission\newline
Bolshoy Karetny 19, Moscow 101447, GSP-4, Russia\newline
olsh@iitp.ru,\  olsh@glasnet.ru

\newpage






\section{Introduction}\label{sec:int}
\setcounter{equation}{0}

The
{\it centralizer construction\/} proposed in \cite{o:ea} shows that 
certain ``quantum"
algebras can be obtained as projective limits
of centralizers in classical enveloping algebras.
This approach has been applied to
the series of matrix Lie algebras 
to construct
the quantum algebras called
the {\it Yangians\/} and {\it twisted Yangians\/}
which are originally defined as certain deformations of
enveloping algebras; cf. \cite{d:ha}.
The type $A$ case is treated in  \cite{o:ea}, \cite{o:ri}
and extended to the
$B,\ts C,\ts D$ types
in \cite{o:ty} and \cite {mo:cc}. A modified version of the $A$ type
construction is given in \cite{m:yt}.

In the $A$ type case, one considers the centralizer
$\A_m(n)$ of the subalgebra $\gl(n-m)$ in
the enveloping algebra $\U(\gl(n))$. It turns out that
for each pair $n>m$ there is a natural algebra homomorphism 
$\A_m(n)\to \A_m(n-1)$
so that one can define a projective limit algebra $\A_m$ by using
the chain of homomorphisms
\beq\label{chaingln}
\cdots\longrightarrow \A_m(n){\longrightarrow} 
\A_m(n-1)\longrightarrow\cdots 
\longrightarrow \A_m(m+1)
{\longrightarrow} \A_m(m).
\end{equation}
The algebra $\A_m$ has a large center $\A_0$ which is isomorphic
to the algebra of shifted symmetric functions $\Lambda^*$ (see \cite{oo:ss})
and one has an isomorphism
\beq\label{isomyang}
\A_m\simeq \A_0\ot \Y_m,
\end{equation}
where $\Y_m$ is the Yangian for the Lie algebra $\gl(m)$;
\cite[Theorem~2.1.15]{o:ri}.
In particular, for each $n\geq m$ there is a natural
homomorphism
\beq\label{homimy}
\Y_m\to \A_m(n).
\end{equation}
This result was used in \cite{c:ni, nt:ry, m:yt} to study
the class of representations of $\Y_m$ which naturally
arises from this construction. 
A similar result for the
$B,\ts C,\ts D$ types \cite{o:ty, mo:cc}
was used in \cite{m:br, m:wb, m:wbgt} to construct weight bases
of Gelfand--Tsetlin type for representations of
the classical Lie algebras.

Our aim in the present paper is to extend
these constructions to the series of the
symmetric groups $S(n)$. Denote by $B_m(n)$ the centralizer
of the subgroup $S(n-m)$ in the group algebra $\C[S(n)]$, where
$S(n-m)$ consists of the permutations which fix
each of the indices $1,2,\dots,m$. It turns out
that no natural analog of the chain \eqref{chaingln} exists
for the algebras $B_m(n)$. Indeed, note that $B_n(n)=B_{n-1}(n)=\C[S(n)]$
and so, by analogy with \eqref{chaingln}
we would have a homomorphism  $\C[S(n)]\to \C[S(n-1)]$ identical
on   $\C[S(n-1)]$. But such a homomorphism does not
exist for $n>4$.

On the other hand, it was observed in \cite{o:ea}
that an analog of \eqref{homimy} still exists:
for any $n\geq m$ there is a homomorphism
$\HH_m\to B_m(n)$, where $\HH_m$
is the {\it degenerate affine 
Hecke algebra\/} introduced by
Drinfeld \cite{d:da} and Lusztig \cite{l:ah}.
This fact was used by Okounkov and Vershik \cite{ov:na}
to develop a new approach to the representation theory of the symmetric
groups; see also earlier results by Cherednik \cite{c:sb}.

This observation together with the semigroup method
\cite{o:ur} allows one to expect that
an analog of the centralizer construction for
the symmetric group should exist, with the
algebras $\C[S(n)]$ replaced with the {\it semigroup algebras\/}
$A(n)=\C[\Gamma(n)]$, where $\Gamma(n)$ is the {\it semigroup
of partial bijections\/} of the set $\{1,\dots,n\}$.
Alternatively, the elements of $\Gamma(n)$ can be identified with
$(0,1)$-matrices which have at most one $1$ in each row and column.
The semigroups $\Gamma(n)$ are studied in \cite{o:ur} and used to prove
Lieberman's classification theorem \cite{l:sc} for unitary
representations of the complete infinite symmetric group.

Taking the centralizer $A_m(n)$ of $\Gamma(n-m)$ in $A(n)$
instead of the algebras $B_m(n)$ we do obtain an analog
of the chain of homomorphisms \eqref{chaingln}
for the algebras $A_m(n)$. The corresponding
projective limit algebra $A_m$ has a decomposition analogous
to \eqref{isomyang}
\beq\label{isomhec}
A_m\simeq A_0\ot \wH_m,
\end{equation}
where $\wH_m$ is a `semigroup analog' of the degenerate	affine 
Hecke algebra $\HH_m$. 
The algebra $\wH_m$ can be presented by generators and
defining relations. Moreover, the algebra $\HH_m$
is a homomorphic image of $\wH_m$. The subalgebra $A_0\subseteq A_m$
is commutative and isomorphic to the algebra of
shifted symmetric functions $\Lambda^*$; see
\cite{oo:ss} for a detailed study of the algebra $\Lambda^*$.

The mentioned above homomorphisms $\HH_m\to B_m(n)$ can be regarded
as `retractions' of the homomorphisms $\wH_m\to A_m(n)$
whose existence is provided with the centralizer
construction for the symmetric groups.

Finally, the algebra $A$ is defined as the inductive limit
of $A_m$ as $m\to \infty$. We show that $A$ naturally acts
in the so-called {\it tame\/} representations of $S(\infty)$
and it can be regarded as the `true' analog
of the group algebra  $\C[S(\infty)]$. Indeed, contrary to the
finite case of $\C[S(n)]$, the algebra $\C[S(\infty)]$
has a trivial center while $A$ contains a large
center whose elements act by scalar operators
in the tame representations of $S(\infty)$.
Moreover, the central elements separate 
the irreducible tame representations.

Some other generalizations of the degenerate affine 
Hecke algebra $\HH_m$ have been constructed by Nazarov
\cite{n:yo, n:ys}.

The paper is organized as follows. 
Section~\ref{sec:pre} is preliminary. Here we
define tame representations of  $S(\infty)$
and describe the properties of the semigroups $\Gamma(n)$;
most of these results are contained in \cite{o:ur}.
In Section~\ref{sec:cc}
we construct the algebras $A_m$ as projective limits
of the centralizers $A_m(n)$. Section~\ref{sec:sa}
describes the algebra $A_0$ and establishes its isomorphism
with the algebra of shifted symmetric
functions $\Lambda^*$. The main results are given
in Section~\ref{sec:sam} where we investigate
the structure of $A_m$
and describe its relationship with
the degenerate affine Hecke algebras.

The authors acknowledge the financial support of
the Australian Research Council. G.O. was also supported
by the Russian Foundation for Basic Research under
grant 98-01-00303.

\section{Tame representations and semigroups}\label{sec:pre}
\setcounter{equation}{0}

Here we give constructions
of the {\it tame\/} representations of
the infinite symmetric group and describe
the semigroups of partial bijections. 
The material of this section is based on 
\cite{o:ur}	and \cite[\S2]{o:nl}.
About applications of the semigroup method to representations of
``big'' groups, see Olshanski \cite{o:sr}, Neretin \cite{n:cs}.

\subsection{Constructing tame representations}\label{subsec:ctr}

Let $\N=\{1,2,\ldots\}$ and $n\in \N$. 
We denote by $S(n)$ the group of 
permutations of the 
set $\N_n:=\{1,\ldots,n\}$. 
We also regard $S(n)$
as the group of permutations of $\N$ fixing $n+1,n+2,\ldots$, and set
$$
S(\infty)=\bigcup_{n\geq 1}S(n).
$$
This is the group of all {\it finite} permutations of the set $\N$.

For any $m\leq n$, we denote by $S_m(n)$ the subgroup of permutations in
$S(n)$ fixing $1,\ldots,m$, and set
$$
S_m(\infty)=\bigcup_{n\geq m} S_m(n).
$$
Note that the subgroups $S(n)$ and $S_n(\infty)$
of $S(\infty)$ commute. 
This simple observation will play an important role later.

For a (unitary) representation $T$, we denote by $H(T)$ its (Hilbert) space.

Let $T$ be a unitary representation of the group $S(\infty)$. 
For $n=0,1,\ldots$,
denote by $H_n(T)$ the subspace of $S_n(\infty)$-invariant vectors in $H(T)$.
Since $\{S_n(\infty)\}$ is an descending chain of subgroups, $\{H_n(T)\}$ is
an ascending chain of subspaces. Let
$$
H_{\infty}(T)=\bigcup_{n\geq 0} H_n(T).
$$
Since $S(n)$ and $S_n(\infty)$ commute, the subspace $H_n(T)$ is
invariant with respect to $S(n)$, and so, $H_{\infty}(T)$ 
is invariant with respect to the whole group $S(\infty)$.

\bde\label{def:tame}{\rm%
A unitary representation $T$ of the group 
$S(\infty)$ is said to
be {\it tame} if $H_{\infty}(T)$ is dense in $H(T)$.
}%
\endproof
\ede

Note that for an irreducible representation $T$, this is equivalent to saying
that $H_{\infty}(T)$ is nonzero.

Clearly, the trivial representation of $S(\infty)$ is tame, and another
one--dimensional representation, $s\mapsto {\rm{sgn}}\ts s $, is not
tame. Less trivial examples follow.

\bex{\rm
(i) Let $H$ be the space $l_2$ with its canonical basis
$e_1,e_2,\ldots$ and let $S(\infty)$ operate in $H$ 
by permuting the basis vectors.
The representation $H$ is tame. It is
irreducible and for any $n$ the  subspace $H_n$ is spanned by
$e_1, \ldots,e_n$.
\par
(ii)
The right (or left) regular representation of the group
$S(\infty)$ in the
Hilbert space $l_2\big(S(\infty)\big)$ is not tame.
}%
\endproof
\eex

For any $n=0,1,2,\ldots$ and any partition $\lambda\vdash n$, we will
construct a tame representation $T_{\lambda}$.
First, if $n=0$ then $\lambda=\emptyset$ and $T_{\emptyset}$ is
the one-dimensional trivial representation. Let now $n\geq 1$
and let $\pi_{\lambda}$ denote the irreducible representation of
$S(n)$ corresponding to $\lambda$. Then $T_{\lambda}$ is defined as
the induced representation
\beq\label{tlam}
T_{\lambda}={\rm{Ind}}^{S(\infty)}_{S(n)\times S_n(\infty)}
(\pi_{\lambda}\otimes 1)
\end{equation}
where 1 stands for the trivial representation of $S_n(\infty)$.
The representation $T_{\lambda}$ can be realized as follows. 
Let $\Omega(n)$ denote
the set of injective maps
$
\omega: \N_n\to \N.
$
Define a right action of $S(\infty)$ on $\Omega(n)$ by
\beq\label{sact}
\omega\cdot s=s^{-1}\circ\omega,\qquad 
s\in S(\infty),
\end{equation}
and a left action of $S(n)$ by
\beq\label{tact}
t\cdot\omega=\omega\circ t^{-1},\qquad t\in S(n).
\end{equation}
Note that these two actions commute. Consider the Hilbert space
$l_2\big(\Omega(n),H(\pi_{\lambda})\big)$ of $H(\pi_{\lambda})$-valued
square-integrable functions on $\Omega(n)$, and let $H(n,\lambda)$ be its
subspace formed by the functions $f(\omega)$ such that
\beq\label{func}
f(t\cdot\omega)=\pi_{\lambda}(t)f(\omega), \qquad t\in S(n),
\quad \omega\in\Omega(n).
\end{equation}
The action of $T_{\lambda}$ in $H(n,\lambda)$ is given by
\beq\label{ts}
(T_{\lambda}(s)f)(\omega)=f(\omega\cdot s).
\end{equation}
The space $H(T_{\lambda})$ may now be identified with
$H(n,\lambda)$.
For any $l\geq n$ set
\beq\label{omegal}
\Omega(n,l)=\{\omega\in\Omega(n)\ | \ \omega(\N_n)
\subseteq \N_l\}
\end{equation}
and note that
\beq\label{omcup}
\Omega(n)=\bigcup_{l\geq n}\Omega(n,l)
\end{equation}
Also, set
\beq\label{omegaln}
H'_l(T_{\lambda})=\{f\in H(n,\lambda)\ |\ \supp f 
\subseteq\Omega(n,l)\},
\end{equation}
where $\supp f=\{\omega\in\Omega(n)\ |\ f(\omega)\ne 0\}$. 

\bpr\label{prop:hl} For any $n\in\N$ we have
\beq\label{hl}
H_l(T_{\lambda})=\begin{cases} \{0\} &\text{if $\; l< n$}, \\
H'_l(T_{\lambda}) \qquad&\text{if $\; l\geq n$}.\end{cases}
\end{equation}
\epr

\Proof Since $S_l(\infty)$ acts trivially on $\Omega(n,l)$, we have
$H'_l(T_{\lambda})\subseteq H_l(T_{\lambda})$. Conversely, let
$f\in H_l(T_{\lambda})$. Then the function $\Vert f(\omega)\Vert^2$ 
is constant
on any orbit of the subgroup $S_l(\infty)$ in $\Omega(n)$. Since 
the sum of the
$\Vert f(\omega)\Vert^2$ taken over $\omega\in\Omega(n)$ must be finite, 
we have
$\Vert f(\omega)\Vert=0$ unless the orbit containing $\omega$ is finite. But if
$\omega\notin\Omega(n,l)$, then its orbit is clearly infinite, so that
$f(\omega)=0$. This proves the opposite inclusion
$H_l(T_{\lambda})\subseteq H'_l(T_{\lambda})$. 
\endproof

\bpr\label{prop:ttame} For any $n\in\N$ and any $\lambda\vdash n$,
the representation $T_{\lambda}$ of $S(\infty)$ is tame and irreducible.
\epr

\Proof By Proposition~\ref{prop:hl},
\beq\label{hinf}
H_{\infty}(T_{\lambda})=\bigcup_l H_l(T_{\lambda})=
\bigcup_l H'_l(T_{\lambda}).
\end{equation}
This is the subspace of finitely supported functions in
$H(T_{\lambda})=H(n,\lambda)$ which is clearly dense. So,
$T_{\lambda}$ is tame.

The subspace
$H_n(T_{\lambda})=H'_n(T_{\lambda})$ is both cyclic in $H(T_{\lambda})$ and
irreducible under the action of the subgroup $S(n)$. This implies
that $T_{\lambda}$ is irreducible.
\endproof

We shall identify any partition $\lambda$ with its Young diagram; see
e.g. \cite{m:sf}.
We write $|\lambda|=n$ if $\lambda$ has $n$ boxes.
Given two diagrams $\lambda$ and $\mu$
the notation $\mu\nearrow\lambda$ will mean that $\mu$ can be obtained from
$\lambda$ by removing one box, i.e. $\mu\subset\lambda$ and
$|\mu|=|\lambda|-1$. An {\it infinite
tableau\/} is defined as an infinite chain of diagrams
\beq\label{ytab}
\tau=(\lambda^{(1)}\nearrow\lambda^{(2)}\nearrow\ldots),\qquad
|\lambda^{(n)}|=n.
\end{equation}
Two infinite tableaux will be called {\it equivalent\/} 
if the corresponding
chains of diagrams differ in a finite number of places only.

Given an infinite tableau $\tau$, we may construct 
an inductive limit unitary representation
$\Pi(\tau)$ of the group $S(\infty)$ as follows. By the branching 
rule for the symmetric groups (see e.g. \cite{jk:rt}, \cite{ov:na}),
for any $n$ the representation $\pi_{\lambda^{(n)}}$ occurs in the 
decomposition of  $\pi_{\lambda^{(n+1)}}\downarrow S(n)$ with multiplicity one. 
Hence there
is an infinite chain of embeddings
\beq\label{emb}
\pi_{\lambda^{(1)}}\hra\pi_{\lambda^{(2)}}\hra\ldots
\end{equation}
which are defined uniquely up to scalar multiples. So we may
set
\beq\label{pitau}
\Pi(\tau)={\rm{lim\,ind\,}} \pi_{\lambda^{(n)}}, \qquad n\to\infty.
\end{equation}
One can show that any $\Pi(\tau)$ is irreducible
(cf. \cite[Theorem~2.1]{o:urid} and \cite[\S2.7]{o:nl}), and that
$\Pi(\tau)$ and $\Pi(\tau')$ are
isomorphic if and only if  $\tau$ and $\tau'$ are equivalent.

This construction provides a large class of 
pairwise non-equivalent irreducible
representations of the group $S(\infty)$. We will be only interested in 
some special
representations of this class.
Let $\lambda=(\lambda_1,\dots,\lambda_l)$ 
be an arbitrary diagram. Consider
an infinite tableau $\tau=(\lambda^{(i)})$ such that
\beq\label{tabinf}
\lambda^{(i)}=(i-|\lambda|,\lambda_1,\ldots,\lambda_l) \quad \text{for}
\quad i\geq|\lambda| +\lambda_1,
\end{equation}
and set
$
\Pi_{\lambda}=\Pi(\tau).
$
The representation $\Pi_{\lambda}$ is well defined
since the equivalence class of $\Pi_{\lambda}$ does not depend on the
choice of  $\lambda^{(i)}$ for small values of $i$.

\bpr\label{prop:equiv} 
The representations
$T_{\lambda}$ and $\Pi_{\lambda}$ are isomorphic
for any $\lambda$.
\epr

\Proof For any $l\geq n+\lambda_1$, the natural representation of
$S(l)$ in the space $H_l(T_{\lambda})$ is isomorphic to the induced
representation
\beq\label{indl}
{\rm{Ind}}^{S(l)}_{S(n)\times S(l-n)}(\pi_{\lambda}\otimes 1) 
\end{equation}
where $S(l-n)$ is identified with $S_n(l)$, and $1$ stands for the trivial
representation of $S(l-n)$. This follows immediately from \eqref{hl}.
It is well known \cite{jk:rt} that the representation \eqref{indl} 
is multiplicity free and that its
spectrum  consists of the representations $\pi_{\mu}$ such that 
$\mu\vdash l$ and
\beq\label{strip}
\mu_{i+1}\leq\lambda_i\leq\mu_{i}, \qquad i\geq 1. 
\end{equation}
It follows 
from \eqref{strip} that
$\pi_{\lambda^{(l)}}$ occurs in the decomposition of \eqref{indl}. 
Let
$H^0_l(T_{\lambda})$ denote the corresponding subspace of $H_l(T_{\lambda})$.
It remains to prove that $H^0_m(T_{\lambda})$ is contained in
$H^0_l(T_{\lambda})$ for any $m<l$ provided that $m$ is large enough. 
This follows from the fact that
\beq\label{indm}
\pi_{\lambda^{(m)}}\subset\pi_{\lambda^{(l)}}|^{}_{S(m)}, 
\end{equation}
and
$
\pi_{\lambda^{(m)}}\not\subset \pi_{\mu}|^{}_{S(m)} 
$
if $\mu$ satisfies \eqref{strip} and
$\mu\neq \lambda^{(l)}$.
Indeed, 
property \eqref{indm} follows from the definition of the diagrams
$\lambda^{(l)}$ for large $l$ and the branching rule. Finally, note
that there exists $i$ such that $\lambda_i>\mu_{i+1}$ (otherwise
$\mu=\lambda^{(l)}$). Applying the branching rule again
we complete the proof.
\endproof

\subsection{The semigroup method}\label{subsec:sm}

\bde\label{def:prodg}{\rm Let $X$ be a set.
\par
(i) A {\it partial bijection\/} of $X$ is a bijection $\gamma: D \to R$
between two (possibly empty) subsets of $X$. The subset $D\subseteq X$ is
called the {\it domain\/} of $\gamma$ and denoted by $\dom\gamma$. The
subset $R\subseteq X$ is called the {\it range\/} of $\gamma$ and denoted by
$\range\gamma$. If $x\in X$ belongs to $\dom\gamma$, then we will say that
$\gamma$ {\it is defined on\/} $x$. The set of partial bijections of $X$ is
denoted by $\Gamma(X)$.
\par
(ii) Given $\gamma\in \Gamma(X)$, we define $\gamma^*\in \Gamma(X)$ as the
inverse bijection $\gamma^{-1}: \range\gamma\to \dom\gamma$, so that
$\dom\gamma^*=\range\gamma$ and $\range\gamma^*=\dom\gamma$. The
mapping $\gamma\to \gamma^*$ is involutive: $(\gamma^*)^*=\gamma$.
\par
(iii) Given $\gamma_1, \gamma_2\in\Gamma(X)$ with
$\dom\gamma_i=D_i$ and $\range\gamma_i=R_i$, $i=1,2$,
their {\it product\/} $\gamma_1\gamma_2$ is a partial bijection on $X$
with $D=\dom\gamma_1\gamma_2=\gamma_2^{-1}(D_1\cap R_2)$ and
$R=\range\gamma_1\gamma_2=\gamma_1(D_1\cap R_2)$:
\beq\label{prod}
\gamma_1\gamma_2=(\gamma_1\big|^{}_{\gamma_2(D)} )
\circ(\gamma_2
\big|^{}_{D}).
\end{equation}
That is, $\gamma_1\gamma_2$ is defined on $x\in X$ if and only if
$\gamma_2$ is defined on $x$ and $\gamma_1$ is defined on $\gamma_2(x)$;
then $(\gamma_1\gamma_2)(x)=\gamma_1(\gamma_2(x))$.
}\endproof
\ede

Any $\gamma\in \Gamma(X)$ may be regarded as a {\it relation\/} on $X$. 
The product defined above is the {\it product of relations}; see \cite{cp:at}.
With this product
$\Gamma(X)$ becomes a semigroup, 
called the {\it semigroup of partial bijections\/}
of $X$.

The involution $\gamma\to\gamma^*$ is an anti-homomorphism of $\Gamma(X)$,
so that $\Gamma(X)$ is an {\it involutive semigroup\/}.

The semigroup $\Gamma(X)$ possesses the {\it unity\/} 1, which is the
identity bijection $X\to X$, and the {\it zero\/} 0, which is the (unique)
bijection of the empty subset of $X$ onto itself. Note that
\beq\label{onezer}
1\ts\gamma=\gamma\ts 1=\gamma, \quad 0\ts\gamma=\gamma\ts 0=0, 
\qquad \text{for any}\ 
\gamma\in\Gamma(X).
\end{equation}
For any subset $Y\subseteq X$, let $1_Y\in \Gamma(X)$ denote the identity
bijection $Y\to Y$. In particular, $1_X=1$ and $1_{\emptyset}=0$. Then
$1_Y$ is a self-adjoint idempotent, i.e.,
\beq\label{adj}
(1_Y)^*=1_Y, \quad (1_Y)^2=1_Y.
\end{equation}
Moreover, all idempotents of this type are pairwise commuting.
For any $\gamma\in\Gamma(X)$ we obviously have
\beq\label{inver}
\gamma^*\ts\gamma=1_{\dom\gamma}, \quad \gamma\ts\gamma^*=1_{\range\gamma}.
\end{equation}
The subset of those $\gamma\in\Gamma(X)$ for which
$\dom\gamma=\range\gamma=X$ is clearly identified with the group
$S(X)$ of all permutations of the set $X$.

\medskip
\noindent {\it Remark.} The semigroup $\Gamma(X)$ is a model example
of an {\it inverse semigroup\/} (each of its elements
has an inverse); see \cite{cp:at}. 
The class of inverse semigroups is
very closed to that of groups, and the role of the semigroups of partial
bijections $\Gamma(X)$ is quite similar to that of 
the symmetric groups $S(X)$.
In particular, any inverse semigroup is isomorphic to an involutive
subsemigroup of some $\Gamma(X)$: this is an analog of Cayley theorem; see
\cite{cp:at}. \endproof

There is a convenient realization of partial bijections by
(0,1)-matrices, i.e., by the matrices whose entries are 0 or 1. 
A (0,1)-matrix is said to be {\it monomial\/} if any of its rows or columns
contains at most one 1. Given a set $X$, we will consider (0,1)-matrices
whose rows and columns are labeled by the points of $X$. Then to any
$\gamma\in\Gamma(X)$, we assign a monomial (0,1)-matrix $[\gamma_{xy}]$ as
follows: for $x,y\in X$
\beq\label{gamxy}
\gamma_{xy}=\begin{cases} 1\qquad&\text{if $y\in\dom\gamma$ and $\gamma(y)=x$},\\
0 &\text{otherwise.}\end{cases}
\end{equation}
In particular, the unity $1\in\Gamma(X)$ and the zero $0\in\Gamma(X)$
are represented by the unit and the zero matrices, respectively.
We thus obtain an isomorphism between the semigroup $\Gamma(X)$ and the
semigroup of monomial matrices over $X$ with the usual matrix
multiplication. Note that in this matrix realization, the involution
$\gamma\mapsto\gamma^*$ is represented by the matrix transposition.

\bde\label{def:maps}{\rm
For $Y\subseteq X$ we define three mappings
as follows:
\par
(i) The {\it canonical projection\/} $\theta:\Gamma(X)\to\Gamma(Y)$.
Let $\gamma\in\Gamma(X)$. Then $\theta(\gamma)$ is defined at $y\in Y$ 
if $\gamma$ is defined at $y$ and $\gamma(y)\in Y$. The image of $y$ 
with respect to $\theta(\gamma)$ is
$\gamma(y)$. In matrix terms: the matrix of $\theta(\gamma)$ is
obtained from that of $\gamma$ by striking the rows and columns
corresponding to points of $X\setminus Y$.
\par
(ii) The {\it canonical embedding\/} 
$\phi: \Gamma(Y)\hra \Gamma(X)$. Let
$\gamma\in\Gamma(Y)$. Then $\phi(\gamma)$ is defined at $x\in X$ if
and only if $x\in\dom\gamma\subseteq Y$ or $x\in X\setminus Y$. In the
former case $\phi(\gamma)$ sends $x$ to $\gamma(x)$, and in the latter
case it fixes $x$. In matrix terms, for $x,y\in X$,
\beq\label{phi}
\phi(\gamma)_{xy}=\begin{cases} \gamma_{xy} \qquad&\text{if $x,y\in Y$},\\
\delta_{xy} &\text{otherwise.} \end{cases}
\end{equation}
\par
(iii) The {\it $0$-embedding} $\psi: \Gamma(Y)\hra \Gamma(X)$. Let
$\gamma\in\Gamma(Y)$. Then $\psi(\gamma)$ is obtained by regarding
$\dom\gamma$ and $\range\gamma$ as subsets of $X$. In matrix terms, for
$x,y\in X$,
\beq\label{psi}
\psi(\gamma)_{xy} =\begin{cases}\gamma_{xy}\qquad&\text{if $x,y\in Y$,}\\
0&\text{otherwise.}\end{cases}
\end{equation}
}\endproof
\ede

For a positive integer $n$ we shall denote
by $\Gamma(n)$ the semigroup $\Gamma(\N_n)$.
Given a tame representation $T$ of $S(\infty)$ we now construct
a representation
$\T_n$ of the semigroup $\Gamma(n)$. 
Regarding $S(\infty)$ as the group of infinite monomial matrices
which have only a finite number of $1$'s off the diagonal,
define the map $\theta^{(n)}:S(\infty)\to \Gamma(n)$
which assigns to each infinite matrix its upper left
$n\times n$ submatrix. It is easy to check that
the map $\theta^{(n)}$ is surjective and it induces a bijection
between the set of double cosets 
$S_n(\infty)\backslash S(\infty)/S_n(\infty)$ and $\Gamma(n)$; 
see \cite{o:ur}.

Let $T$ be a tame representation of $S(\infty)$.
Denote by $P_n$ the orthogonal projection $P_n:H(T)\to H_n(T)$.
Suppose that $n$ is so large that $H_n(T)\ne \{0\}$.
For any $\sigma_1,\sigma_2\in S(\infty)$ we have 
\beq\label{repgam}
\theta^{(n)}(\sigma_1)=\theta^{(n)}(\sigma_2)\quad\Rightarrow\quad
P_nT(\sigma_1)P_n=P_nT(\sigma_2)P_n.
\end{equation}
Therefore there exists a unique map $\T_n:\Gamma(n)\to \End H_n(T)$
such that
\beq\label{mapstn}
\T_n\big(\theta^{(n)}(\sigma)\big)=	P_nT(\sigma)|^{}_{H_n(T)},\qquad 
\sigma\in S(\infty).
\end{equation}

\bpr\label{prop:trep}
The map $\T_n=\T_n(T)$ defined by \eqref{mapstn} is 
a representation of the semigroup $\Gamma(n)$
in $H_n(T)$.
\epr

\Proof We give a sketch of the proof. The details and one more proof
can be found in \cite{o:ur}.
We show first that the tame representation $T$ of $S(\infty)$
can be extended to a representation $\T$ of the semigroup
$\Gamma(\infty)$ of partial bijections of the set $\N$.

Further, we prove that $P_n$ coincides with the operator
$\T(1_n)$ where $1_n$ is the identity bijection
of the subset $\N_n$.

Finally, consider the $0$-embedding $\psi: \Gamma(n)\to\Gamma(\infty)$;
see \eqref{psi}. We have for any $\gamma\in\Gamma(n)$
\beq\label{mulipt}
\T_n(\gamma)=P_n\T\big(\psi(\gamma)\big)P_n=
\T(1_n)\T\big(\psi(\gamma)\big)\T(1_n)
=\T(1_n	\psi(\gamma)1_n)=\T(\psi(\gamma)).
\end{equation}
This proves the multiplicativity of $\T_n$.
\endproof

Consider the tame representations $T_{\lambda}$
of $S(\infty)$ constructed
in the previous section.
Proposition~\ref{prop:trep} yields a representation
$\T_n(\lambda):=\T_n(T_{\lambda})$ of
$\Gamma(n)$ provided that $H_n(T_{\lambda})\ne 0$, i.e.,
$n\ge |\lambda|$.

We now outline the
proof of the classification theorem for
representations of $\Gamma(n)$; see \cite{o:ur}.

\bth\label{thm:clasgam}
The representations $\T_n(\lambda)$ where $\lambda$ is
a partition with $|\lambda|\leq n$
exhaust all irreducible representations of $\Gamma(n)$.
\eth

\Proof
Let $\T$ be a representation of $\Gamma(n)$.
Then for any $m\leq n$ the operator $\T(1_m)$
is a projection.  Denote its image by $H_m(\T)$.
We let $m(\T)$ denote the minimum value of $m$ such that
$H_m(\T)\ne\{0\}$.

Further, if $\T$ is irreducible then one shows that $\dim\T<\infty$.
If $m=m(\T)$ then the subspace $H_m(\T)$ is invariant
under $S(m)$ and irreducible. So, $H_m(\T)$ corresponds to
a partition $\lambda$ with $|\lambda|=m$. 
A standard argument proves that $\T$ is uniquely determined by
$\lambda$.

Conversely, given a partition $\lambda$ with $|\lambda|=m\leq n$
one uses the following argument
to construct an irreducible representation $\T$ of $\Gamma(n)$
such that $m=m(\T)$ and the representation of $S(m)$ in
$H_m(\T)$ corresponds to $\lambda$.

Denote by $\Omega(m,n)$ the set of injective maps
$\omega:\N_m\to\N_n$. We define $H(\T)$ to be the space
of functions $f: \Omega(m,n)\to H(\pi_{\lambda})$ such that
\beq\label{repta}
f(t\cdot\omega)=\pi_{\lambda}(t)f(\omega),\qquad t\in S(m),\qquad
\omega\in \Omega(m,n);
\end{equation}
see \eqref{tact}.
The action of $\gamma\in\Gamma(n)$ is given by
\beq\label{gamnact}
(\T(\gamma)f)(\omega)=\begin{cases}
f(\gamma^*\omega)\qquad
&\text{if $\omega=(\omega_1,\dots,\omega_m)
\subseteq \dom\gamma^*$},\\
0\qquad&\text{otherwise}.
\end{cases}
\end{equation}
One easily checks that $\T$ is a representation of $\Gamma(n)$.
Moreover, it is isomorphic to
the representation $\T_n(\lambda)$.
\endproof

Note that the representation of $\Gamma(n)$ corresponding to $m=0$
and the empty diagram is the trivial representation sending all
elements of $\Gamma(n)$ to $1$.

Theorem~\ref{thm:clasgam} leads to the following result.

\bth\label{thm:clirtam}
Let $\lambda$ range over
the set of all Young diagrams
including the empty diagram. The representations $T_{\lambda}$ 
constructed
in Section~\ref{subsec:ctr} exaust, within
equivalence, all the irreducible tame representations of the group
$S(\infty)$. Moreover, any tame representation of $S(\infty)$ can be
decomposed into a direct sum of irreducible tame representations.
\eth

\Proof See Theorem~6.7 and \S7.2 in \cite{o:ur}. \endproof

This is equivalent to Lieberman's theorem \cite{l:sc} concerning
continuous unitary representations of the complete infinite symmetric group
(this group consists of all permutations of the set $\N$);
see \cite[\S7]{o:ur}.

\section{Centralizer construction}
\label{sec:cc}
\setcounter{equation}{0}

We shall denote by $\theta_{n}$ the canonical 
projection
$\Gamma(n)\to\Gamma(n-1)$;
see Definition~\ref{def:maps}.
So, if $\gamma\in \Gamma(n)$
then $\theta_{n}(\gamma)$
is the upper left corner of $\gamma$ of order
$n-1$. Here the elements of $\Gamma(n)$ and $\Gamma(n-1)$ are regarded as
$(0,1)$-matrices of order $n$ and $n-1$, respectively.

We set
$A(n)=\C[\Gamma(n)]$, the semigroup algebra of $\Gamma(n)$. 
The canonical embedding $\Gamma(n-1)\hra\Gamma(n)$
is extended to an embedding $A(n-1)\hra A(n)$ by linearity.
Further, set
$A(\infty)=\C[\Gamma(\infty)]=\bigcup_{n\ge 1}A(n)$, the semigroup algebra of
$\Gamma(\infty)$.

For each $i=1,\dots,n$ denote by $\ve_i$ the diagonal 
$n\times n$-matrix whose
$ii$-th entry is $0$ and all other diagonal entries are equal to $1$.
The corresponding element of $\Gamma(n)$
is the identity bijection of the set $\N_n\setminus \{i\}$.
The algebra  $A(n)$ is obviously
generated by $S(n)$ and the elements $\ve_i$,
$i=1,\dots,n$. We have for any $i$ and any $\sigma\in S(n)$:
\beq\label{epssi}
\sigma\ts\ve_i\ts \sigma^{-1}=\ve_{\sigma(i)}.
\end{equation}

\bpr\label{prop:defream}
The algebra $A(n)$ is isomorphic
to the abstract algebra with generators $s_1,\dots,s_{n-1}$,
$\ve_1,\dots,\ve_n$ and the defining relations
\begin{align}\label{defsn}
s_k^2&=1,\qquad &s_k\ts s_{k+1}\ts s_k&=s_{k+1}\ts s_k\ts s_{k+1},\qquad
s_k\ts s_l=s_l\ts s_k,\quad |k-l|>1,\\
\label{defga}
\ve_k^2&=\ve_k,\qquad &\ve_k\ts\ve_l&=\ve_l\ts\ve_k,\\
\label{defsga}
s_k\ts\ve_k&=\ve_{k+1}\ts s_k,\qquad &s_k\ts \ve_k\ts\ve_{k+1}&=\ve_k\ts\ve_{k+1}.
\end{align}
\epr

\Proof Denote the abstract algebra by $\A(n)$.
The assignments $s_k\mapsto (k,k+1)$ and $\ve_k\mapsto\ve_k$
obviously define an algebra epimorphism $\A(n)\to A(n)$.
Note also that \eqref{defsn} are defining relations for the symmetric
group $S(n)$ and \eqref{epssi} holds. 
To complete the proof we verify that
$\dim \A(n)\leq \dim A(n)$. We have
\beq\label{diman}
\dim A(n)=|\ts\Gamma(n)|=\sum_{r=0}^n\binom{n}{r}^2 r!\ .
\end{equation}
On the other hand, we see from the relations 
\eqref{defga} and \eqref{defsga} that
\beq\label{decamc}
\A(n)=\C[S(n)]\ts\C[\ve_1,\dots,\ve_n].
\end{equation}
Here $\C[\ve_1,\dots,\ve_n]$ is the subalgebra of $\A(n)$
generated by the $\ve_k$. It is spanned by the
monomials $\ve_{k_1}\cdots \ve_{k_r}$ with $k_1<\cdots<k_r$.
Given such a monomial,
consider the subspace $\C[S(n)]\ts\ve_{k_1}\cdots \ve_{k_r}$
in $\A(n)$. Using \eqref{epssi} if necessary,
we may assume without loss of generality that
$k_i=i$ for each $i$. Observe that by \eqref{defsga} we have in $\A(n)$
\beq\label{coset}
\sigma\ts S(r)\ts\ve_1\cdots \ve_r = \sigma\ts\ve_1\cdots \ve_r,
\qquad \sigma\in S(n).
\end{equation}
Hence the dimension of the subspace $\C[S(n)]\ts\ve_1\cdots \ve_r$
does not exceed
the number of left cosets of $S(n)$ over the subgroup $S(r)$.
Therefore, $\dim\A(n)$ does not exceed
\beq\label{dimacn}
\sum_{r=0}^n \binom{n}{r}\frac{n!}{r!}
\end{equation}
which coincides with \eqref{diman}.
\endproof

Using Proposition~\ref{prop:defream}
we shall sometimes identify $A(n)$ with the algebra $\A(n)$.

\bco\label{cor:retra}
The mapping 
\beq\label{retr}
(k,k+1)\mapsto (k,k+1),\qquad \ve_k\mapsto 0
\end{equation}
defines an algebra homomorphism $A(n)\to \C[S(n)]$ which is
identical on the subalgebra $\C[S(n)]$. \endproof
\eco

We shall call \eqref{retr} the {\it retraction homomorphism\/}.
It can be equivalently defined as follows.
For any $\gamma\in\Gamma(n)$, define its {\it rank\/}, 
denoted by $\rank\gamma$, as the number of $1$'s in the
$(0,1)$-matrix representing $\gamma$. That is,
\beq\label{rank}
\rank\gamma= 
|\ts\dom\gamma|=|\ts\range\gamma|.
\end{equation} 
The rank of an element $a=\sum a_{\gamma}\gamma\in A(n)$
is defined as the maximum of the ranks $\rank\gamma$ 
with $a_{\gamma}\ne 0$.
Now \eqref{retr} can also be defined by
setting for $\gamma\in \Gamma(n)$
\beq\label{retrdef}
\gamma\to \begin{cases} \gamma\qquad&\text{if}\ \rank\gamma=n,\\
0\qquad&\text{if}\ \rank\gamma<n,
\end{cases}
\end{equation}
and extending this to $A(n)$ by linearity. 

For any $0\leq m\leq n$ 
denote by $\Gamma_m(n)$ the subsemigroup of $\Gamma(n)$
which consists of the matrices with first $m$ diagonal entries
equal to $1$. 
Set
$A_m(n)=A(n)^{\Gamma_m(n)}$, the centralizer of
$\Gamma_m(n)$ 
in the algebra $A(n)$. 
In particular, $A_0(n)$ is the center of $A(n)$.

We extend
$\theta_{n}$ to a linear map $A(n)\to A(n-1)$.

\bpr\label{prop:hom} 
The restriction of $\theta_{n}$
to $A_{n-1}(n)\subseteq A(n)$
defines a unital algebra homomorphism
\beq\label{hom}
\theta_{n}: A_{n-1}(n)\to A(n-1).
\end{equation}
Moreover,
\beq\label{homam}
\theta_{n}(A_{m}(n))\subseteq A_{m}(n-1)\qquad 
\text{for\ \  $m=0,1,\ldots,n-1$.}
\end{equation}
\epr

\Proof Note that any $a\in A_{n-1}(n)$
commutes with $\ve_n$
because $\ve_n$ is contained in $\Gamma_{n-1}(n)$, and 
$A_{n-1}(n)$ is the centralizer of $\Gamma_{n-1}(n)$ in $A(n)$.
Since $\ve_n$ is an idempotent we have
\beq\label{pra}
\theta_n(a)=a\ts\ve_n=\ve_n\ts a=\ve_n\ts a\ts\ve_n,\qquad 
a\in A_{n-1}(n).
\end{equation}
Now, if $a,b\in A_{n-1}(n)$ then
\beq\label{thab}
\theta_n(ab)= \ve_n\ts ab\ts\ve_n =	\ve_n\ts a\ts \ve_n\ts\ve_n\ts b\ts\ve_n
=\theta_n(a)\theta_n(b).
\end{equation}

It is clear that \eqref{hom} preserves the unity.

Finally, we need to show that if $a\in A_m(n)$ and 
$b\in \Gamma_{m}(n-1)$ then
$\theta_{n}(a)$ and $b$ commute. Regard $b$ as an element of $\Gamma(n)$; 
then it
lies in $\Gamma_m(n)$ and commutes with $\ve_n$. This implies that 
$a\ts\ve_n$
and $b$ commute, and so do $\theta_{n}(a)$ and $b$.
\endproof

For $\gamma\in\Gamma(n)$, set
\beq\label{jgamma}
\bsp
J(\gamma)&=\{i\ |\ 1\leq i\leq n,\  \gamma_{ii}=0\}, \\ 
\deg\gamma&=
|\ts J(\gamma)|.
\esp
\end{equation}

\bpr\label{prop:degga} Let $\gamma,\gamma' \in\Gamma(n)$. Then
\beq\label{gamgam}
\deg\gamma\gamma'\leq\deg\gamma+\deg\gamma'.
\end{equation}
Moreover, the equality in \eqref{gamgam} implies 
$\gamma\gamma'=\gamma'\gamma$.
\epr

\Proof Regard $\gamma$ and $\gamma'$ as partial bijections of the set
$\N_n=\{1,\ldots,n\}$. If $\delta\in \Gamma(n)$ then
$\N_n\setminus J(\delta)$ is the set of
$\delta$-invariant elements in the domain of $\delta$. 
This implies
\beq\label{caps}
\big(\N_n\setminus J(\gamma)\big)\cap
\big(\N_n\setminus J(\gamma')\big)\subseteq\big(\N_n\setminus 
J(\gamma\gamma')\big).
\end{equation}
Therefore,
\beq\label{capuni}
J(\gamma)\cup J(\gamma')\supseteq J(\gamma\gamma'),
\end{equation}
and \eqref{gamgam} follows.
Finally, the equality in \eqref{gamgam} implies that $J(\gamma)$ and
$J(\gamma')$ are disjoint. But then $\gamma$ and $\gamma'$ must commute.
\endproof

Using \eqref{jgamma}, we define a filtration of the space $A(n)$:
\beq\label{filtr}
\C=A^0(n)\subseteq A^1(n)\subseteq\ldots \subseteq A^n(n)=A(n)
\end{equation}
where $A^M(n)$ is spanned by
the subset $\{\gamma\ |\  \deg\gamma\leq M\}\subseteq \Gamma(n)$. By
Proposition~\ref{prop:degga} 
this filtration is compatible with the algebra structure of
$A(n)$, and the corresponding graded algebra 
\beq\label{grad}
\gr A(n)=\bigoplus_M\ts\big(A^M(n)/A^{M-1}(n)\big)
\end{equation}
is commutative.
Note that for $\gamma\in\Gamma(n)$ the degree
$
\deg\theta_{n}(\gamma)
$
can be equal either to
$\deg\gamma$ or  $\deg\gamma-1$. Therefore
the homomorphisms \eqref{homam} are compatible with the
filtration on $A(n)$.

\bde\label{def:algam}{\rm
For $m=0,1,2,\ldots$ let $A_m$ be the projective limit of the
infinite sequence 
\beq\label{seqal}
\cdots\longrightarrow A_m(n)\overset{\theta_n}{\longrightarrow} 
A_m(n-1)\longrightarrow\cdots 
\longrightarrow A_m(m+1)
\overset{\theta_{m+1}}{\longrightarrow} A_m(m)
\end{equation}
taken in the category of filtered algebras.
}\endproof
\ede

By the definition, an element $a\in A_m$ is a sequence $(a_n\ |\ n\geq m)$ 
such that
\beq\label{seqa}
a_n\in A_m(n),\quad\theta_{n}(a_n)=a_{n-1}, \qquad\deg a :=\sup_{n\geq
m}\deg a_n <\infty
\end{equation}
with the componentwise operations. For $n\geq m$ we shall  denote by
$\theta^{(n)}$ the projection $A_m\to A_m(n)$ such that
\beq\label{thetn}
\theta^{(n)}(a)=a_n.
\end{equation}
The $M$-th term of the filtered algebra $A_m$ will be denoted
by $A^M_m$. 

There are natural algebra homomorphisms 
$A_m\to A_{m+1}$ defined by
\beq\label{amin}
(a_n\ |\ n\geq m) \mapsto (a_n\ |\ n\geq m+1)
\end{equation}
where we use the inclusions $A_m(n)\subset A_{m+1}(n)$ for $n>m$. 
These homomorphisms are injective because $a_m$ 
is uniquely determined by
$a_{m+1}$.

\bde\label{def:indlim} 
{\rm
The algebra $A$ is defined as the inductive limit (the union) of
the algebras $A_m$ taken with respect to the embeddings $A_m\hra
A_{m+1}$, $m\geq 0$, defined in \eqref{amin}.
}
\endproof
\ede

Since these embeddings preserve the filtration, $A$ is a filtered algebra. We
will denote by $A^M$ the $M$-th term of the filtration, so that
\beq\label{mfilt}
A^M=\bigcup_{m\geq 0} A^M_m.
\end{equation}

\bpr\label{prop:stabemb} There exists a natural embedding
$A(\infty)\hra A$ whose image consists of stable
sequences $a=(a_n)\in A$.
\epr

\Proof Let $b\in A(\infty)$. There exists $m$ such that $b\in
A(m)$. Note that $b\in A_m(n)$ for any $n\geq m$ since $A(m)$ and
$\Gamma_m(n)$ commute. Set
\beq\label{abseq}
a=(a_n\ |\  n\geq m)\in A_m\subset A\quad\text{with $a_n\equiv b$.}
\end{equation}
The sequence $a\in A$ only depends on $b$ 
and not on the choice of $m$. The mapping
$b\mapsto a$ is clearly an algebra embedding.
\endproof

\bco\label{cor:embsinf} There is a natural algebra embedding
$\C[S(\infty)]\hra A$. \endproof
\eco

\bpr\label{prop:center} The center of the algebra $A$ coincides with
$A_0$.
\epr

\Proof Recall that $A_0(n)$ is the center of $A(n)$. 
The subalgebra $A_0$ is contained in the center of $A$
since the sequences $a=(a_n)\in A$ are multiplied componentwise.

Conversely, if $a$ belongs to the center of $A$ then $a$ 
commutes with the subalgebra
$A(\infty)\subset A$. This implies that for any $n$
the element $a_n$ is contained in
$A_0(n)$, and so $a\in A_0$.
\endproof

\medskip
\noindent {\it Remark.} The same argument
shows that the subalgebra
$A_m\subset A$ coincides with the centralizer in $A$ of the subalgebra
\beq\label{subgm}
\bigcup_{n\geq m}\C\lbrack\Gamma_m(n)\rbrack\subset A(\infty)\subset A.
\end{equation}
\endproof

Note that the centers of both algebras
$\C[S(\infty)]$ and $A(\infty)$ are trivial.
However, as it will be shown in the next section, the center $A_0$ of the
algebra $A$ has a rich structure.

\bpr\label{prop:amod} 
For any tame representation $T$ of the group
$S(\infty)$, the subspace $H_{\infty}(T)\subseteq H(T)$ admits a natural structure
of an $A$-module such that for any $m$ the subspace $H_m(T)$ is invariant
with respect to the subalgebra $A_m$ 
{\rm (}and hence $H_n(T)$ is invariant
with respect to $A_m$ for $n\geq m$\/{\rm )}.
\epr

\Proof Let $a\in A$ and $h\in H_{\infty}(T)$. Choose $m$ such that
$a\in A_m$. Then we may write $a=(a_n\ |\  n\geq m)$. Let
us prove that
\beq\label{stable}
a_nh=a_{n+1}h,\quad n\geq m.
\end{equation}
Consider the family of representations
$\{\T_n\}$ associated with $T$, which has been
introduced in Section~\ref{subsec:sm}. 
Each $\T_n$ is a representation of the semigroup
$\Gamma(n)$ in the space $H_n(T)$ and so it 
can be extended to a representation of the
semigroup algebra $A(n)$ in the same space. Recall that
$\T_{n+1}(\ve_{n+1})$ projects  $H_{n+1}(T)$ 
onto its subspace $H_n(T)$.
Since $h$ is already contained in $H_n(T)$ (as we assume $n\geq m$), we have
$\T_{n+1}(\ve_{n+1})h=h$, so that
\beq\label{teps}
T_{n+1}(1-\ve_{n+1})h=0.
\end{equation}
This implies that $h$ is annihilated by the left ideal $I(n+1)\subset A(n+1)$.
Since $a_n-a_{n+1}\in I(n+1)$, this implies \eqref{stable}.

Define a mapping
\beq\label{ahtoh}
A\times H_{\infty}(T)\to H_{\infty}(T), \qquad (a,h)\mapsto a_mh 
\end{equation}
where $m$ is so large that 
$a\in A_m$ and $h\in H_m(T)$. Note that
under this assumption  $ah\in H_m(T)$.

The mapping \eqref{ahtoh} is clearly bilinear and $1\ts h=h$. 
The
multiplicativity property $(ab)h=a(bh)$ follows from
the definition of the multiplication in $A_m$. 
\endproof

\bpr\label{prop:aoscal} If{\ }  $T$ is an irreducible tame representation of
$S(\infty)$ then $H_{\infty}(T)$ is irreducible as an $A$-module. 
In particular, the center $A_0$ of $A$ acts by scalar operators.
\epr

\Proof The first claim is obvious because $H_{\infty}(T)$ is already
irreducible as a $A(\infty)$-module (see Proposition~\ref{prop:stabemb}).
To prove the second claim consider an element $a\in A_0$ as an operator
in $H_{\infty}(T)$.
It suffices to show that $a$ has an eigenvalue.	The result will then follow
by a standard argument using Schur's lemma.

Assume that $a$ has no eigenvalues. Let us show first that
$a$ is algebraically independent over $\C$. Indeed,
let $P(x)\in\C[x]$ be a nonzero polynomial
of a minimum degree such that $P(a)=0$.
Then $P(x)=(x-\alpha)\ts Q(x)$ for a certain $\alpha\in\C$ and
a polynomial $Q(x)\in\C[x]$. Since $Q(a)$ is a nonzero
operator, there is a vector $v\in H_{\infty}(T)$ such that
$w:=Q(a)v\ne 0$. Then $w$ is an eigenvector for $a$ with the eigenvalue
$\alpha$. Contradiction.

We note now that the space $H_{\infty}(T)$ has countable dimension
and then use a version of Dixmier's argument \cite{d:ri} as follows.

Since $a$ is algebraically independent over $\C$, the field
$\C(a)$ is embedded in the endomorphism algebra of $H_{\infty}(T)$.
This implies that the dimension of $H_{\infty}(T)$ is at least as large
as the dimension of $\C(a)$ over $\C$, but the latter
is continuum. This contradiction completes the proof.
\endproof

\section{The structure of the algebra $A_0$}
\label{sec:sa}
\setcounter{equation}{0}

In the last two sections we aim to describe the structure
of the algebras $A_m$.
Here we consider the commutative
algebra $A_0$; see Proposition~\ref{prop:center}.
We construct generators of $A_0$ and show that
it is isomorphic to the algebra of shifted symmetric
functions.

\subsection{Generators of $A_0$}
\label{subsec:genao}

Let $Z(S(n))$ denote the center of the algebra $\C\lbrack S(n)\rbrack$.
For $0\leq M\leq n$ denote by 
$Z^{M}(S(n))$ the $M$-th term of the filtration	on $Z(S(n))$
inherited from the algebra $A(n)$; see \eqref{filtr}.
Note that
\beq\label{zozi}
Z^0(S(n))=Z^1(S(n))=\C\ts 1
\end{equation}
because, for a permutation $s\in S(n)$
the inequality $\deg s\leq 1$ implies $s=1$.

For any partition 
$\M=(M_1,\ldots,M_r)$
with 
\beq\label{modm}
|\M|=M_1+\ldots +M_r\leq n 
\end{equation}
introduce the element $c_n^{\M}$ of the group algebra $\C\lbrack S(n)\rbrack$
as follows
\beq\label{cnm}
c_n^{\M}=
\sum(i_1,\ldots,i_{M_1})(j_1,\ldots,j_{M_2})
\cdots  (k_1,\ldots,k_{M_r}) 
\end{equation}
where the sum is taken over the sequences
$i_1,\ldots,i_{M_1};j_1,\ldots,j_{M_2}; \ldots; 
k_1,\ldots,k_{M_r}$ of $|\M|$ pairwise distinct indices taken from
$\N_n$.	By $(i_1,\ldots,i_{M_1})$ etc. in \eqref{cnm}
we denote cycles in the symmetric group. 
For the empty partition $\emptyset$ we set $c_n^{\emptyset}=1$. 
Note that 
$
c_n^{(1)}=\sum_{i=1}^n (i)=n\ts 1.
$
Given two partitions $\M$ and $\LL$ we denote by $\M\cup\LL$
the partition whose parts are those of $\M$ and $\LL$ 
rewritten in the decreasing order.
We have
\beq\label{cnmones}
c_n^{\M\cup 1\cup\ldots\cup 1}=
(n-\vert\M\vert)\cdots(n-\vert\M\vert-p+1)\ts c_n^{\M}
\end{equation}
where $p$ stands for the number of $1$'s in the left hand side of 
the relation.

By definition \eqref{jgamma} of the degree of an element
of $\Gamma(n)$ we have
\beq\label{degrecn}
\deg c_n^{(1)}=0,\qquad\text{and}\qquad\deg c_n^{(M)}=M\quad
\text{for $M\geq 2$.}
\end{equation}
More generally,
\beq\label{degcnm}
\deg c_n^{(M_1,\ldots,M_r)}=\sum_{i,\; M_i\geq 2}M_i.
\end{equation}

\bpr\label{prop:bascent} Each of the families
\beq\label{bas1}
c_n^{\M},\qquad\vert\M\vert=n, 
\end{equation}	
and  
\beq\label{bas2}
c_n^{\M}, \qquad
\vert\M\vert\leq n\ \  \text{and $\M$ has no part equal to
$1$,}
\end{equation}
forms a basis of
$Z(S(n))$. Moreover, the elements of degree $\leq M$ of each family
form a basis of $Z^{M}(S(n))$. 
\epr

\Proof The elements \eqref{bas1} are proportional to the characteristic
functions of the conjugacy classes of the 
group $S(n)$ and so, they form a basis of $Z(S(n))$. 
By \eqref{cnmones} the elements of type \eqref{bas2} 
are proportional to those of type
\eqref{bas1}.
\endproof

\bpr\label{prop:prodcn} 
For any two partitions
$\M=(M_1,\ldots,M_r)$ and $\LL=(L_1,\ldots,L_t)$ with
$\vert\M\vert+\vert\LL\vert\leq n$ we have
\beq\label{prodcn}
c_n^{\M}c_n^{\LL}=c_n^{\M\cup\LL}\ts +\ts(\ldots),
\end{equation}
where $(\ldots)$ stands for a linear combination of the elements
$c_n^{{\cal K}}$ with $\vert{\cal K}\vert<\vert\M\vert+\vert\LL\vert$.
\epr

\Proof For a permutation $s\in S(n)$ or $s\in S(\infty)$ 
define its {\it support\/}
as
\beq\label{supp}
\supp s=\{i\in \N_n\ |\  s(i)\ne i\}\quad\text{or} \quad\supp s=\{i\in
\N\ |\  s(i)\ne i\},
\end{equation}
respectively.
(The degree of a permutation is then given by
$\deg s =|\ts\supp s|$; cf. \eqref{jgamma}).

Let $s\in S(n)$ be a permutation which occurs
in the expansion of $c_n^{\M}$, that is, $s$ is of cycle type
$\M\cup 1\cup\ldots\cup 1$ (with $n-|\M|$ units). 
Similarly, let $s'$ be a permutation
occurring in $c_n^{\LL}$. If the supports
$\supp s$ and $\supp s'$ are disjoint
then $s$ and $s'$ commute, and the product $ss'$ occurs in the
expansion of $c_n^{\M\cup\LL}$. In particular,
$\deg ss'=\vert\M\vert+\vert\LL\vert$.
If $\supp s$ and $\supp s'$ have a non-empty
intersection then the degree of
$ss'$ is strictly less than $\vert\M\vert+\vert\LL\vert$.
\endproof

\noindent
{\it Remark.} A detailed investigation
of the structure constants for the
products of type \eqref{prodcn} have been recently
given by Ivanov and Kerov \cite{ik:ac}.
\endproof

\bco\label{cor:gencent} 
Let $k=(k_1,\ldots,k_n)$ run over the $n$-tuples of 
non-negative integers such that 
$2k_2+\cdots +nk_n\leq n$.
Then the monomials
\beq\label{gencent}
(c_n^{(2)})^{k_2}\cdots (c_n^{(n)})^{k_n}
\end{equation}
form a basis of $Z(S(n))$. Moreover, for any $M\geq 0$, the elements
\eqref{gencent} 
with $2k_2+\cdots +nk_n\leq M$ form a basis of $Z^{M}(S(n))$.
\eco

\Proof It suffices to prove that
\beq\label{ckcm}
(c_n^{(2)})^{k_2}\cdots (c_n^{(n)})^{k_n}=c_n^{\M}\ts +\ts (\ldots)
\end{equation}
where
$
\M=2^{k_2}\cdots n^{k_n}
$
and $(\ldots)$ stands for a certain linear combination of the elements 
$c_n^{\M'}$ with $\vert\M'\vert< \vert\M\vert=2k_2+\cdots +nk_n$. 
But this
follows from Proposition~\ref{prop:prodcn}.
\endproof

Now we will define analogs of the elements $c_n^{\M}$ for the
algebra $A_0(n)$. Namely, for any partition $\M=(M_1,\ldots,M_r)$ with
$\vert\M\vert\leq n$ set
\beq\label{deltanm}
\Delta_n^{\M}=\sum(i_1,\ldots,i_{M_1})(j_1,\ldots,j_{M_2}) \cdots 
(k_1,\ldots,k_{M_r}) (1-\epsi_{i_1})\cdots(1-\epsi_{k_{M_r}} )
\end{equation}
where, as in \eqref{cnm}, the sum is taken over all sequences
of $|\M|$ pairwise distinct indices taken from $\N_n$.
In particular, 
\beq\label{delta1}
\Delta_n^{(1)}=\sum_{i=1}^n (1-\epsi_i).
\end{equation}
For the empty partition $\emptyset$	 we set 
$\Delta_n^{\emptyset}=1$.
By \eqref{jgamma}, we have
\beq\label{degdel}
\deg\Delta_n^{\M}=\vert\M\vert\qquad \text{for any partition $\M$},
\end{equation}
cf. \eqref{degcnm}.
Note that $\Delta_n^{\M}$ can also be written as
\beq\label{del1}
\Delta_n^{\M}=\sum (1-\epsi_{i_1})\cdots(1-\epsi_{k_{M_r}})
(i_1,\ldots,i_{M_1})(j_1,\ldots,j_{M_2}) \cdots 
(k_1,\ldots,k_{M_r}),
\end{equation}
and as
\beq\label{del2}
\bsp
\Delta_n^{\M}=\sum (1-\epsi_{i_1})\cdots(1-\epsi_{k_{M_r}})
&(i_1,\ldots,i_{M_1})(j_1,\ldots,j_{M_2}) \\
{}\cdots{} 
&(k_1,\ldots,k_{M_r})(1-\epsi_{i_1})
\cdots(1-\epsi_{k_{M_r}}).
\esp
\end{equation}
Indeed, 
$(1-\epsi_{i_1})\cdots(1-\epsi_{k_{M_r}})$
is invariant under the conjugation by the permutation 
$(i_1,\ldots,i_{M_1})(j_1,\ldots,j_{M_2}) \cdots (k_1,\ldots,k_{M_r})$
which implies \eqref{del1}. To derive \eqref{del2} it suffices to note that 
$(1-\epsi_{i_1})\cdots(1-\epsi_{k_r})$ is an idempotent.

\bpr\label{prop:delao} 
The element $\Delta_n^{\M}$ belongs to $A_0(n)$ for any $\M$.
\epr

\Proof 
Since $\Gamma(n)$ is generated by the group $S(n)$ and the pairwise commuting
idempotents $\epsi_1,\ldots,\epsi_n$, it suffices to show that
$\Delta_n^{\M}$ commutes both with $S(n)$ and with the $\epsi_i$. 
The first
claim is clear since $\Delta_n^{\M}$ is invariant under the conjugation by 
the elements of $S(n)$. 
To prove the second claim, we observe that any $\epsi_l$,
$1\leq l\leq n$, commutes with any term 
\beq\label{term}
\sigma=(i_1,\ldots,i_{M_1})(j_1,\ldots,j_{M_2}) \cdots 
(k_1,\ldots,k_{M_r}) (1-\epsi_{i_1})\cdots(1-\epsi^{}_{k_{M_r}})
\end{equation}
in \eqref{deltanm}. Indeed, this is clear if $l$ does not occur in the set
of indices in \eqref{term} because 
$\epsi_l$ commutes with the corresponding
cycle. But if $l$ coincides with one of the indices $i_1,\ldots,k_{M_r}$, 
then $\epsi_l\ts \sigma=\sigma\ts \epsi_l=0$. 
This follows from \eqref{del2} and
the relation $(1-\epsi_l)\epsi_l=0$. 
\endproof

\bpr\label{prop:thetdel} We have
\beq\label{thetdel}
\theta_{n}(\Delta_n^{\M})=\Delta_{n-1}^{\M}
\end{equation}
where we adopt the convention that 
\beq\label{delmo}
\Delta_k^{\M}=0\quad\text{if} \quad \vert\M\vert>k.
\end{equation}
\epr

\Proof By the definition of the projection
$\theta_{n}$ (see Section~\ref{sec:cc})
we need to calculate $\Delta_n^{\M}\ts\epsi_n$. 
However, as it follows from the proof of Proposition~\ref{prop:delao}, 
the effect of multiplying 
$\Delta_n^{\M}$ by $\epsi_n$ reduces to striking from \eqref{deltanm} 
all terms \eqref{term}
such that $n$ occurs among the
 corresponding indices. If $\vert\M\vert=n$, then
all the terms are vanished, so that the result of the multiplication 
is $0$. If
$\vert\M\vert<n$, then the terms that survive are just the terms of the sum
defining $\Delta_{n-1}^{\M}$. 
\endproof

\bco\label{cor:deltam} For any partition $\M$, there exists an element
$\Delta^{\M}\in A_0$  such that
\beq\label{deltam}
\theta^{(n)}(\Delta^{\M})=\Delta_n^{\M}\quad \text{for any $n\geq 1$}
\end{equation}
with the convention \eqref{delmo}.
\eco

\Proof By Proposition~\ref{prop:delao}, $\Delta_n^{\M}\in A_0(n)$. Now we
apply Proposition~\ref{prop:thetdel} and note 
that the degrees of the elements $\Delta_n^{\M}$ are uniformly bounded by
\eqref{degdel}. 
\endproof

We now aim to prove an analog of Proposition~\ref{prop:bascent}
for the algebra $A_0(n)$; see  Proposition~\ref{prop:basao}
below. For this we need the following three lemmas.

Let $I(n)=A(n)(1-\ve_n)$ denote the left ideal of the algebra $A(n)$
generated by the element $1-\ve_n$.

\ble\label{lem:inao}  For any $n$,
\beq\label{inao}
I(n)\cap A_0(n)=Z(S(n))(1-\epsi_1)\cdots(1-\epsi_n). 
\end{equation}
\ele

\Proof First suppose that $x\in A(n)$ can be written as 
$y\ts (1-\epsi_1)\cdots(1-\epsi_n)$ where $y\in Z(S(n))$.
The argument of the proof of Proposition~\ref{prop:delao}
shows that $x\in A_0(n)$.
Moreover, we obviously have $x\in I(n)$.

Conversely, suppose $x\in I(n)\cap A_0(n)$. Then $x\ts\epsi_n=0$. 
Using the invariance of $x$ under the conjugation by elements of $S(n)$
we also obtain $x\ts\epsi_i=0$ for $i=1,\ldots,n$. Therefore $x$ is 
invariant under the right
multiplication by $(1-\epsi_1)\cdots(1-\epsi_n)$.
Further, we may write
\beq\label{xyeps}
x=y+\sum_{r=1}^n\sum_{1\leq i_1<\cdots<i_r\leq n}y_{i_1\ldots i_r}
\ts\epsi_{i_1}\cdots\epsi_{i_r},
\end{equation}
where $y$ and all the $y_{i_1\ldots i_r}$ are elements
of $\C\lbrack S(n)\rbrack$; see Proposition~\ref{prop:defream}.
Multiplying this relation by $(1-\epsi_1)\cdots(1-\epsi_n)$ on the right
we obtain
\beq\label{xy1eps}
x=y\ts (1-\epsi_1)\cdots(1-\epsi_n). 
\end{equation}
Finally, for any $s\in S(n)$ we may write
\beq\label{sxs}
x=sxs^{-1}=sys^{-1}(1-\epsi_1)\cdots(1-\epsi_n).
\end{equation}
Averaging over $s\in S(n)$ turns $y$ into an element of $Z(S(n))$.
\endproof

For a subset $I=\{i_1,\dots,i_k\}$ in $\N_n$ we put
$\ve^{}_I=\ve^{}_{i_1}\cdots \ve^{}_{i_k}$, and for
$s\in S(n)$  set
\beq\label{qs}
Q(s)=\{i\in \N_n\ |\ s_{ii}=1\}=\N_n\setminus \supp s. 
\end{equation}

\ble\label{lem:sngn} The mapping
\beq\label{sngn}
s\mapsto \gamma,\qquad \gamma= s\ts \epsi^{}_{Q(s)}=\epsi^{}_{Q(s)}s, 
\end{equation}
defines a bijection of $S(n)$ onto the set of all $\gamma\in\Gamma(n)$
satisfying the conditions
\begin{align}\label{dr}
\dom\gamma&=\range\gamma,  \\
\label{degn}
\deg\gamma&=n. 
\end{align}
\ele

\Proof The effect of the multiplication of $s$ by $\epsi^{}_{Q(s)}$  
from the left or from the right consists of replacing all the 1's on the
diagonal by zeros. This implies \eqref{dr}, and \eqref{degn} is obvious.

Conversely, let $\gamma\in\Gamma(n)$ satisfy \eqref{dr} and \eqref{degn}. 
Relation \eqref{dr} means that
for any $i\in \N_n$ the $i$-th
row and the $i$-th column are zero or non-zero at the same time, whereas
\eqref{degn} means that all the diagonal entries of $\gamma$ are zero.
Now, let the matrix $\sigma$ be defined as follows.  
Set $\sigma_{ii}=1$ if  the $i$-th row (and
the $i$-th column) of $\gamma$ is zero, and set $\sigma_{ij}=\gamma_{ij}$ for
$i\ne j$.
It is easy to see that $\sigma\in S(n)$ and that $\gamma$ is the
image of $\sigma$ under the mapping \eqref{sngn}.
\endproof

\ble\label{lem:inj} 
The restriction of the
projection $\theta_{n}: A_0(n)\to A_0(n-1)$ to the subspace $A_0^{n-1}(n)$
is injective.
\ele

\Proof Let $x\in A_0(n)$ and $\theta_{n}(x)=0$. We will show that
$\deg x=n$ unless $x=0$.
By Lemma~\ref{lem:inao}, $x$ can be written 
as a linear combination of the elements
\beq\label{seps}
s\ts (1-\epsi_1)\cdots(1-\epsi_n)=\sum_{I\subseteq\ts \N_n}
(-1)^{|I|} s\ts \epsi^{}_I, \qquad s\in S(n). 
\end{equation}
Rewrite this as
\beq\label{seps2}
s\ts (1-\epsi_1)\cdots(1-\epsi_n)= 
\sum_{I\supseteq Q(s)} (-1)^{|I|}s\ts\epsi^{}_I +\sum_{I\not\supseteq Q(s)} 
(-1)^{|I|} s\ts\epsi^{}_I.
\end{equation}
Then the terms of the first sum are of degree $n$ whereas those of the second 
sum are
of degree  strictly less than $n$. So it suffices to prove that the elements
\beq\label{linind}
\sum_{I\supseteq Q(s)} (-1)^{|I|} s\ts\epsi^{}_I, \qquad s\in S(n), 
\end{equation}
are linearly independent. 
Note that
\begin{align}
\label{rkse1}
\rank s\ts\epsi^{}_I &=n-|Q(s)| \qquad\text{if $I=Q(s)$,}\\
\label{rkse2}
\rank s\ts\epsi^{}_I &<n-|Q(s)| \qquad\text{if $I\supset Q(s)$};
\end{align}
see \eqref{rank}.
Write $S(n)$ as the disjoint union of $n+1$ subsets:
\beq\label{uni}
S(n)=\bigcup_{k=0}^n \{s\in S(n)\ |\ n-|Q(s)|=k\}.
\end{equation}
If $s$ belongs to the $k$-th subset then
\beq\label{rsum}
\rank\left(\sum_{I\supseteq Q(s)} (-1)^{|I|}s\ts\epsi^{}_I\right)=k. 
\end{equation}
Moreover, only one term of the sum in \eqref{rsum} has rank $k$, namely that
with $I=Q(s)$.
Finally, it remains to note that by Lemma~\ref{lem:sngn} all the elements
$s\epsi^{}_{Q(s)}$ with $s\in S(n)$ are pairwise distinct elements of
$\Gamma(n)$. 
\endproof

\bpr\label{prop:basao} For any $n$ the elements 
$\Delta_n^{\M}$, where $\M$ is any partition with $\vert\M\vert\leq n$,
form a basis of $A_0(n)$. Furthermore, for any $M$ such that 
$0\leq M\leq n$ these
elements with $\vert\M\vert\leq M$ form a basis of $A^M_0(n)$.
\epr

\Proof The first claim of the proposition will follow from
the second one. We will prove the second claim using induction on
$n$.
The claim is obviously true for $n=1$. 
Assume that $n\geq 2$ and $M\leq n-1$. By the induction hypothesis the elements
$\Delta_{n-1}^{\M}$ with $\vert\M\vert\leq M$ form a basis of $A_0^{M}(n-1)$.
By Proposition~\ref{prop:thetdel} the image of $\Delta_n^{\M}$ under
$\theta_{n}$ is
$\Delta_{n-1}^{\M}$. By Lemma~\ref{lem:inj} the restriction
$\theta_{n}\downarrow
A_0^{M}(n)$ is injective. Therefore, the elements $\Delta_n^{\M}$
with $\vert\M\vert\leq M$ form a basis in $A_0^{M}(n)$.

Further, let us show that the elements $\Delta_n^{\M}$ 
with $\vert\M\vert=n$ 
form a basis of $I(n)\cap A_0(n)$.
Note that
\beq\label{delcm}
\Delta_n^{\M}=c_n^{\M}\ts(1-\epsi_1)\cdots (1-\epsi_n). 
\end{equation}
By Proposition~\ref{prop:bascent} the elements $c_n^{\M}$, where
$\M$ runs over the set of  partitions of $n$, form a basis of $Z(S(n))$.
Due to
Lemma~\ref{lem:inao}
it now
remains to check that the elements $c_n^{\M}$, being multiplied by
$(1-\epsi_1)\cdots (1-\epsi_n)$, remain linearly independent. 
However, this follows
from the fact that the composite map
\beq\label{compos}
\C[S(n)]\to A(n)\to \C[S(n)]
\end{equation}
is the identity map;
here the first arrow is the multiplication by $(1-\epsi_1)\cdots
(1-\epsi_n)$, and the second arrow is the retraction 
homomorphism \eqref{retr}.

Finally, let us show that
\beq\label{aodec}
A_0(n)=A_0^{n-1}(n)\oplus \big(I(n)\cap A_0(n)\big). 
\end{equation}
Indeed,	as it was shown above, $\theta_{n}$ maps 
$A_0^{n-1}(n)$ onto $A_0^{n-1}(n-1)=A_0(n-1)$ . Since $I(n)\cap A_0(n)$
is the kernel of the restriction $\theta_{n}\downarrow A_0(n)$ and since
$\theta_{n}(A_0(n))$ is contained in $A_0(n-1)$, we obtain
the decomposition
\beq\label{aos}
A_0(n)=A_0^{n-1}(n)+ \big(I(n)\cap A(n)\big). 
\end{equation}
Lemma~\ref{lem:inj} implies that
\beq\label{inters}
A_0^{n-1}(n)\cap I(n)=\{0\} 
\end{equation}
and \eqref{aodec} follows.

To complete the proof we need to show
that the elements $\Delta_n^{\M}$ with $0\leq |\M|\leq n$
form a basis of $A_0(n)$. However, 
the elements with $\vert\M\vert<n$ form a basis of the first
component of the decomposition \eqref{aodec}, 
whereas the
elements with $\vert\M\vert=n$ form a basis in the second component of this
decomposition. 
\endproof

The following is an analog of Corollary~\ref{cor:gencent}.

\bco\label{cor:genaon} 
Let $k=(k_1,\ldots,k_n)$ run over the $n$-tuples of
non-negative integers such that 
$
k_1+2 k_2+\cdots+nk_n\leq n.
$
Then the monomials
\beq\label{delmon}
(\Delta_1^{(1)})^{k_1}\cdots(\Delta_n^{(n)})^{k_n} 
\end{equation}
form a basis of $A_0(n)$. Moreover, for any $M\geq 0$, the monomials 
\eqref{delmon} with
$k_1+2 k_2+\cdots+nk_n\leq M$ form a basis of $A_0^M(n)$.
\eco

\Proof It suffices to prove that
\beq\label{delpar}
(\Delta_1^{(1)})^{k_1}\cdots(\Delta_n^{(n)})^{k_n}=\Delta_n^{\M}+(\ldots)
\end{equation}
where
$
\M=1^{k_1}2^{k_2}\cdots n^{k_n}
$
and $(\ldots)$ stands for a linear combination of the elements 
$\Delta_n^{\M'}$ with $\vert\M'\vert<\vert\M\vert$. Then our
claim will follow from Proposition~\ref{prop:basao}.
To prove \eqref{delpar} we verify that for any 
partitions $\M=(M_1,\ldots,M_r)$ and
$\LL=(L_1,\ldots,L_t)$ with $\vert\M\vert+\vert\LL\vert\leq n$ 
\beq\label{unipart}
\Delta_n^{\M}\Delta_n^{\LL}=\Delta_n^{\M\cup\LL}+(\ldots),
\end{equation}
where the rest term $(\ldots)$ has degree strictly less than
$\vert\M\vert+\vert\LL\vert$ and so, it is a
linear combination of elements  
$\Delta_n^{\K}$ with $\vert\K\vert<\vert\M\vert+\vert\LL\vert$.
Write
\beq\label{deldel}
\Delta_n^{\M}=\sum\delta^{}_I,\qquad \Delta_n^{\LL}=\sum\delta'_J.
\end{equation}
Here $I$ is a sequence $i^{}_1,\ldots,i^{}_{|\M|}$ of pairwise distinct
indices taken from $\N_n$ and
\beq\label{deltai}
\delta^{}_I=(i^{}_1,\ldots,i^{}_{M_1})\ldots (i^{}_{M_1+\ldots+M_{r-1}+1},
\ldots,i^{}_{\vert\M\vert})\ts\prod_{p=1}^{\vert\M\vert}(1-\epsi_{i_p});
\end{equation}
the $\delta'_J$ are the corresponding elements for the partition $\LL$.
Then
\beq\label{demdel}
\Delta_n^{\M}\Delta_n^{\LL}=\sum_{I,\ts J}\delta^{}_I\delta'_J= \sum_{I\cap
J=\emptyset}\delta^{}_I\delta'_J + \sum_{I\cap
J\ne\emptyset}\delta^{}_I\delta'_J. 
\end{equation}
The first sum on the right hand side of \eqref{demdel} is 
$\Delta_n^{\M\cup\LL}$ 
whereas the second sum
is of degree strictly less than $\vert\M\vert+\vert\LL\vert$.
\endproof

Consider the elements $\Delta^{\M}\in A_0$ 
introduced in Corollary~\ref{cor:deltam}.
We shall denote by $\PP$ the set of all partitions.

\bth\label{thm:basao} The elements $\Delta^{\M}$, $\M\in\PP$
form a basis of the algebra $A_0$. Moreover, for any $M\geq 0$,
the elements $\Delta^{\M}$ with $\vert\M\vert\leq M$ form a basis of the
$M$-th subspace $A_0^{M}$ in $A_0$.
\eth

\Proof The first claim follows from the second one. The second
claim follows from Proposition~\ref{prop:basao}
and the definition of $A_0^{M}$ as the projective limit of the
spaces $A_0^{M}(n)$.
\endproof

\bco\label{cor:isomao} For $n>M$, the mapping
\beq\label{isomaon}
\theta_{n}: A_0^M(n)\to A_0^M(n-1)
\end{equation}
is an isomorphism of vector spaces and so is
the mapping
\beq\label{isomao}
\theta^{(n)}: A_0^M\to A_0^M(n), \qquad	 n\geq M.
\end{equation}
In particular, $\dim A_0^M<\infty$.	\endproof
\eco

\bth\label{thm:basano}
The monomials
\beq\label{moninfin}
(\Delta^{(1)})^{k_1}(\Delta^{(2)})^{k_2}\cdots
\end{equation}
with $k_1,k_2, \ldots \in \Z_+$	and  $k_1+2k_2+\ldots <\infty$
form a basis of the algebra $A_0$. Moreover, for any $M\geq 0$
the monomials \eqref{moninfin} with $k_1+2k_2+\ldots \leq M$
form a basis of the
subspace $A_0^M$.
\eth

\Proof It suffices to check that
\beq\label{equivdel}
(\Delta^{(1)})^{k_1}(\Delta^{(2)})^{k_2}\ldots\equiv\Delta^{\M}\quad
\mod A_0^{M-1}
\end{equation}
where $\M=1^{k_1}2^{k_2}\cdots$ and $M=\vert\M\vert$. 
However, this follows from the relation
\eqref{delpar}.
\endproof

\bco\label{cor:geneao} The elements 
$\Delta^{(1)}, \Delta^{(2)}, \ldots$ are
algebraically independent and
generate the algebra $A_0$. \endproof
\eco

\subsection{The algebra $\Lambda^*$ of shifted symmetric functions, 
and the isomorphism
$A_0\simeq \Lambda^*$}\label{subsec:ass}

Let $\Lambda^*(n)\subseteq \C[x_1,\ldots,x_n]$ denote the
subalgebra of
polynomials in $n$ variables $x_1,\ldots,x_n$ which are symmetric in
the new variables
\beq\label{shifts}
y_1=x_1-1,\ \ y_2=x_2-2,\ts \ldots,\ts y_n=x_n-n.
\end{equation}

Following \cite{oo:ss}, we refer to $\Lambda^*(n)$ as the {\it algebra of
shifted symmetric polynomials\/} in $n$ variables. 
We equip $\Lambda^*(n)$ with the
filtration with respect to the usual degree of polynomials.
Set $\Lambda^*(0)=\C$ and for $n\geq 1$ define the projection 
$\Lambda^*(n)\to \Lambda^*(n-1)$ by
specializing $x_n=0$. Note that this projection preserves
the filtration.

\bde\label{def:ss}{\rm  The {\it algebra $\Lambda^*$ of
shifted symmetric functions\/} is the projective limit 
of the {\it filtered}
algebras $\Lambda^*(n)$ as $n\to\infty$.
}\endproof
\ede

In other words, an element $f\in \Lambda^*$ is a 
sequence $(f_n\ |\  n\geq 0)$  such
that
\par
(i) $f_n\in \Lambda^*(n)$ for any $n$; 
\par
(ii) for any $n\geq 1$, $f_n\mapsto f_{n-1}$ under the
projection $\Lambda^*(n)\to \Lambda^*(n-1)$;
\par
(iii) $\deg f_n$ remains bounded as $n\to\infty$.

For an element $f=(f_n)\in \Lambda^*$, we define its {\it degree\/} by
\beq\label{degf}
\deg f=\sup_n\deg f_n,
\end{equation}
and for $M=0,1,\ldots$ we denote by $(\Lambda^*)^M$ 
the subspace in $\Lambda^*$ consisting
of the elements of degree ${}\leq M$. The algebra $\Lambda^*$
was first introduced in \cite{o:urid}. A detailed study
of $\Lambda^*$ is contained in \cite{oo:ss}.

Note an evident similarity between the shifted symmetric functions and
the symmetric functions. Recall (see \cite{m:sf}) that 
the algebra $\Lambda$ of
symmetric functions is defined as the projective limit as $n\to\infty$ of the
{\it graded\/} algebras
$\Lambda(n)\subseteq \C[x_1,\ldots,x_n]$ of symmetric polynomials in $n$
variables. A difference between $\Lambda^*$ and $\Lambda$
consists in a shift of variables and the replacement of the gradation by
a filtration. 
The algebra $\Lambda^*$ may be viewed as a deformation of the algebra
$\Lambda$. Indeed, let $h$ be a numerical parameter, 
and let $\Lambda^*_h$ be defined
similarly to $\Lambda^*$ but with $y_i=x_i-ih$ instead of \eqref{shifts}. 
Then
the algebras $\Lambda^*_h$ with $h\ne 0$, are 
naturally isomorphic to each other.
Moreover,
$\Lambda^*_1$ coincides with $\Lambda^*$ while
$\Lambda^*_0$ coincides with $\Lambda$.
Another relation between $\Lambda^*$ and $\Lambda$ is given by

\bpr\label{prop:gradla} The graded algebra
\beq\label{gradla}
\gr \Lambda^*=\C\oplus\;\bigoplus_{M=1}^{\infty}
\Big((\Lambda^*)^M/(\Lambda^*)^{M-1}\Big)
\end{equation}
is isomorphic to the algebra $\Lambda$.
\epr

\Proof For any $M\geq 1$ and any $n$,
$\Lambda^*(n)^M/\Lambda^*(n)^{M-1}$ is
naturally isomorphic to the $M$-th homogeneous component of the algebra
$\Lambda(n)\subseteq\C[x_1,\ldots,x_n]$. Moreover, this isomorphism is
compatible with the projections 
$\Lambda^*(n)\to \Lambda^*(n-1)$ and $\Lambda(n)\to
\Lambda(n-1)$. This yields an isomorphism $\gr \Lambda^*\to \Lambda$. 
\endproof

In the following example we give some families of generators
of the algebra $\Lambda^*$.  Note that there also exist
other important families analogous  
to the basic symmetric functions; see \cite{oo:ss}.   

\bex\label{def:classym}{\rm
For $M=1,2,\ldots$, 
elements $e_M$, $h_M$, and  $p_M$
defined by the formulas below, are 
shifted symmetric functions:
\begin{align}
E(t)&=1+\sum_{M=1}^{\infty}e_Mt^M=\prod_{k=1}^{\infty}
\frac{1+(x_k-k)t}{1-kt},
\non\\
H(t)&=1+\sum_{M=1}^\infty
h_Mt^M=\prod_{k=1}^{\infty}\frac{1+kt}{1-(x_k-k)t},
\non\\
p_M&=\sum_{k=1}^{\infty}\left((x_k-k)^M-(-k)^M\right).
\non
\end{align}
}\endproof
\eex

The generating functions satisfy the following relations; 
cf. \cite{m:sf}: 
\beq\label{releh}
E(t)H(-t)=1,\qquad \sum_{k=1}^{\infty}p_Mt^M=t\frac{d}{dt}\log H(t).
\end{equation}

\bpr\label{prop:gener} 
The algebra $\Lambda^*$ is isomorphic to the algebra of
polynomials in countably many generators. Furthermore, we have
\beq\label{generla}
\Lambda^*=\C[e_1,e_2,\ldots]=\C[h_1,h_2,\ldots]=\C[p_1,p_2,\ldots].
\end{equation}
\epr

\Proof The corresponding statement
for the algebra $\Lambda$ of symmetric functions
is well known, see \cite[Ch.~1,~Section~2]{m:sf}. Now, we apply
Proposition~\ref{prop:gradla} 
and note that the image of the shifted symmetric function $e_M$,
$h_M$ or $p_M$ in the space 
$(\Lambda^*)^M/(\Lambda^*)^{M-1}\simeq \Lambda^M$
is the corresponding
$M$-th symmetric function (elementary, complete or power
sum). This implies that each of the three families is algebraically
independent and generates the algebra $\Lambda^*$.
\endproof

Let $\Fun \PP $ denote the algebra of complex functions on the
set of partitions $\PP$. By Propositions~\ref{prop:amod}
and \ref{prop:aoscal}, 
there is an algebra homomorphism
\beq\label{pphom}
A_0\to \Fun\PP, \qquad a\mapsto \wh a, 
\end{equation}
such that for $a\in A_0$ and $\lambda\in\PP$, the element $a$ acts in
$H_{\infty}(T_{\lambda})$ as the scalar operator $\wh a(\lambda)\cdot 1$.
On the other hand, any $\lambda\in\PP$ can be viewed as a
sequence $(\lambda_1,\lambda_2,\ldots,0,0,\ldots)$ with finitely many
non-zero coordinates, and so, any element of 
$\Lambda^*$ may be viewed as a function on
$\PP$. Thus we obtain an algebra homomorphism 
$\Lambda^*\to\Fun\PP$ which is clearly an
embedding.

Let $\lambda$ be a partition with $m=|\lambda|\leq n$.
Consider the corresponding
irreducible representation $\pi_{\lambda}$ 
of $S(m)$, and the representation $\T_n(\lambda)$
of the semigroup $\Gamma(n)$; see 
Section~\ref{subsec:sm}.

\bpr\label{prop:eige}
The eigenvalue of the central element $\Delta^{(r)}_n$
in $\T_n(\lambda)$ is $0$ if $r>m$. If $r\leq m$ then
the eigenvalue coincides with that of the element $c^{(r)}_m$
in the representation $\pi_{\lambda}$ of $S(m)$.
\epr

\Proof Recall the construction of $\T_n(\lambda)$
given in Section~\ref{subsec:sm}. Let $\omega$ be an injective map
from $\{1,\dots,m\}$ to $\{1,\dots,n\}$. Regarding $\omega$
as an $m$-tuple $\omega=(\omega_1,\dots,\omega_m)$
we have 
\beq\label{epsaf}
\ve_a\ts f(\omega)=\begin{cases}
0\qquad&\text{if $a\in\omega$},\\
f(\omega)\qquad&\text{if $a\notin\omega$}.
\end{cases}
\end{equation}
Therefore, the product $(1-\ve_{i_1})\cdots (1-\ve_{i_r})$
is a projection to the subspace of functions $f$ such that
the indices $i_1,\dots,i_r$ belong to any $\omega\in \supp f$.
This implies the first statement. The second follows from
the obvious embedding $H(\pi_{\lambda})\subseteq H(\T_n(\lambda))$
whose image consists of the functions supported by the
maps $\omega$ such that 
$\{\omega_1,\dots,\omega_m\}=\{1,\dots,m\}$.
\endproof

It was proved in \cite{ko:pf} (see also \cite{oo:ss})
that the eigenvalue of $c^{(r)}_n$
in the irreducible representation $\pi_{\lambda}$ of $S(n)$
is a shifted symmetric function whose highest homogeneous
component is the power sum symmetric function $p_r$.

\bth\label{thm:alaisom} 
Let $\Lambda^*$ be identified with its image in $\Fun\PP$.
Then the mapping \eqref{pphom} is an isomorphism 
$A_0\to \Lambda^*$ of filtered algebras.
\eth

\Proof By Proposition~\ref{prop:eige}
the images of the generators $\Delta^{(r)}\in A_0$
with respect to the homomorphism \eqref{pphom}
are shifted symmetric functions which are
algebraically independent generators of 
the algebra $\Lambda^*$. The map obviously respects
the filtrations.
\endproof

Recall that by Propositions~\ref{prop:amod} and \ref{prop:aoscal},
elements of the center
$A_0$ act in irreducible tame representations of $S(\infty)$ by scalar
operators. Hence, any such representation determines a homomorphism
$A_0\to \C$.

\bco\label{cor:sep} 
The center $A_0$ separates irreducible tame
representations of $S(\infty)$. That is, non-equivalent
irreducible tame representations give rise to distinct homomorphisms
$A_0\to\C$.
\eco

\Proof  By Theorem~\ref{thm:clirtam}, the irreducible tame representations are
precisely the representations $T_\lambda$. Hence, our claim is
equivalent to the fact that the map $\Lambda^*\to\Fun\PP$ defined
above is an embedding. 
\endproof

\section{The structure of the algebra $A_m$, $m>0$}\label{sec:sam}
\setcounter{equation}{0}

Here we generalize the results of Section~\ref{sec:sa} 
to the algebra $A_m$, where
$m=1,2,\ldots$. Throughout the section we assume $0\leq m\leq n$ and use
the notation
\beq\label{nmn}
\N_{mn}=\{m+1,\ldots,n\}.
\end{equation}
For $\gamma\in\Gamma(n)$, set
\beq\label{jmgam}
\bsp
J_m(\gamma)&=\{i\ |\  i\in\N_{mn},\  \gamma_{ii}=0\},\\
\deg_m\gamma&=|\ts J_m(\gamma)|.
\esp
\end{equation}
We shall call $\deg_m\gamma$ the $m$-{\it degree\/} of $\gamma$.

\bpr\label{prop:degm} For $\gamma,\delta\in\Gamma(n)$,
\beq\label{degm}
\deg_m\gamma\delta\leq\deg_m\gamma+\deg_m\delta.
\end{equation}
\epr

\Proof For any $i\in \N_{mn}$ we have
\beq\label{gamdeli}
(\gamma\delta)_{ii}=0\quad \Rightarrow\quad 
\gamma_{ij}\delta_{ji}=0\qquad\text{for all
$j=1,\ldots,n$.}
\end{equation}
In particular, $(\gamma\delta)_{ii}=0$ implies 
$\gamma_{ii}\delta_{ii}=0$, i.e.,
\beq\label{jmsub}
J_m(\gamma\delta)\subseteq J_m(\gamma)\cup J_m(\delta),
\end{equation}
and \eqref{degm} follows. 
\endproof

\bde\label{def:mfilt}{\rm 
Using the $m$-degree we define a new
filtration in $A(n)$, called the $m$-{\it filtration\/}, by
\beq\label{mffilt}
A(m)=F^0_m(A(n))\subseteq F^1_m(A(n))\subseteq \ldots \subseteq
F^{n-m}_m(A(n))=A(n). 
\end{equation}
Here $F^M_m(A(n))$, the $M$-th term of the filtration, is formed by the
elements $a\in A(n)$ which are 
linear combinations of the
elements of $\Gamma(n)$ of $m$-degree $\leq M$. 
For any subspace $S$ of $A(n)$ we will use
the symbol $F^M_m(S)$ to indicate the $M$-th term of the induced
filtration. 
}\endproof
\ede

By Proposition~\ref{prop:degm}, 
the $m$-filtration is compatible with the algebra
structure of $A(n)$, so the corresponding graded algebra exists.  But
contrary to the case $m=0$, this graded algebra is {\it not\/} commutative
for $m\ge 1$ since it contains, as the 0-component, the
non-commutative algebra $A(m)$.

Let $D$ be a multiplicative
semigroup with unity 1. 
Consider the union $D\cup\{0\}$, where 0 is an extra symbol,
and adopt the convention that
\beq\label{convo}
d\ts 0=0\ts d=0,\quad d+0=0+d=d\quad\text{for any $d\in D$}. 
\end{equation}

\bde\label{def:wrpr}{\rm (i) The semigroup $S(m, D)$ consists of the
$m\times m$ matrices $\alpha=[\alpha_{ij}]$ with entries in $D\cup \{0\}$
such that any row and column contains exactly one non-zero entry.
The product is the matrix multiplication with the conventions
\eqref{convo}.
\par
(ii) The semigroup $\Gamma(n,D)$ is defined as in (i)
by allowing any row and column contain {\it at most\/} one non-zero entry.
}\endproof
\ede

Note that if $D=\{1\}$, then $S(n, D)$ and $\Gamma(n,D)$
coincide with $S(n)$ and $\Gamma(n)$, respectively. If $D$ is a group,
then $S(n,D)$ is the wreath product of $S(n)$ and $D$.

We shall be assuming now that $D$ is the free abelian semigroup 
$\{1,z,z^2,\ldots\}$ with unity 1 and one generator $z$. 
This semigroup is isomorphic to the additive semigroup $\Z_+$. 
We denote the corresponding semigroups introduced in
Definition~\ref{def:wrpr} by
$S(m,\Z_+)$ and $\Gamma(m,\Z_+)$.

Set $\ord z^k=k$ for $k=0,1,\ldots$, and
for $\alpha\in\Gamma(m,\Z_+)$, set
\beq\label{ord}
\ord\alpha=\sum_{i,j;\ts\alpha_{ij}\ne 0} \ord\alpha_{ij}. 
\end{equation}

\bde\label{def:snm}{\rm
(i) Set
\beq\label{snm}
\Gamma(m,n)=\{\sigma\in\Gamma(n)\ |\ \text{$\dom\sigma$ and $\range\sigma$
contain $\N_{mn}$}\}.
\end{equation}
\par
(ii) Consider the linear span of $\Gamma(m,n)$ and
let $Z_m(n)\subset A(n)$ denote the subspace 
in this span formed by the elements invariant under the conjugation
by the elements of the group $S_m(n)$.
}\endproof
\ede

In particular, $\Gamma(0,n)=S(n)$ and $Z_0(n)=Z\big(S(n)\big)$ is the center of
$\C[S(n)]$.
The role of $\Gamma(m,n)$ and
$Z_m(n)$ will be similar to that of $S(n)$ and
$Z\big(S(n)\big)$ in Section~\ref{sec:sa}. 
Note also that $Z_m(n)$ contains $\C[S(n)]^{S_m(n)}$, the centralizer of
$S_m(n)$ in the group algebra $\C[S(n)].$

Now our purpose is to construct a convenient basis in $Z_m(n)$. To do
this, we need to classify the $S_m(n)$-orbits in $\Gamma(m,n)$ where 
the elements of $S_m(n)$ 
act by conjugations.

\bpr\label{prop:param} There is a natural parameterization of the
$S_m(n)$-orbits in $\Gamma(m,n)$ by the couples $(\alpha,\M)$,
where $\alpha\in  \Gamma(m,\Z_+)$ and $\M$ is a partition
such that
\beq\label{param}
\ord\alpha+|\M|=n-m. 
\end{equation}		  
\epr

\Proof Fix an arbitrary element $\sigma\in \Gamma(m,n)$ and assign to it
an $m\times m$-matrix $\alpha=\alpha(\sigma)$ as follows. For $i,j\notin
\N_{mn}$ set
\begin{align}\label{alphaij}
&\text{$\alpha_{ij}=0$} &{}&\text{if $j\notin \dom\sigma$,} \\
&\text{$\alpha_{ij}=1$} &{}&\text{if $j\in\dom\sigma$ and $\sigma(j)=i$,} \\
&\text{$\alpha_{ij}=z^k$} &{}&\text{if $j\in\dom\sigma$}, 
\end{align}
and there exist $k$ points $p_1,\ldots,p_k\in\N_{mn}$ such that
$\sigma(j)=p_1$, $\sigma(p_1)=p_2$, $\ldots$, $\sigma(p_{k-1})=p_k$,
$\sigma(p_k)=i$.
Thus, to any $j\in\dom\sigma$ with $\sigma(j)\in \N_{mn}$ we have
assigned a subset $\{p_1,\ldots,p_k\}\subseteq \N_{mn}$. It is clear that these
subsets are pairwise disjoint. Let $P=P(\sigma)$ denote their union. Then
$
\ord\alpha=|\ts P|\leq n-m.
$
It is also clear that $\alpha\in\Gamma(m,\Z_+)$.

Further, let $P^*=P^*(\sigma)$ be the complement of $P$ in $\N_{mn}$.
Then $P^*$ is contained in the domain of $\sigma$, and 
$P^*$ is $\sigma$-invariant. Therefore,
the restriction of $\sigma$ to $P^*$ defines a permutation of $P^*$. Let
$\M=\M(\sigma)$
be the partition of the number
$|\ts P^*|$ which is defined by the lengths of the
cycles of this permutation. Then the couple
\beq\label{couple}
(\alpha,\M)=(\alpha(\sigma),\M(\sigma)) 
\end{equation}
satisfies \eqref{param}.
It is clear that the couple \eqref{couple} 
remains unchanged if $\sigma$
is replaced by $s\sigma s^{-1}$ with $s\in S_m(n)$. 
Moreover, it is also
clear that if the couples \eqref{couple} 
corresponding to two elements of $\Gamma(m,n)$ are
the same, then these elements belong to the
same orbit. Finally, any couple
satisfying \eqref{param} can be obtained 
from an element of $\Gamma(m,n)$.
\endproof

\noindent
{\it Remark.} A couple \eqref{couple} corresponds to an element of
$S(n)\subseteq \Gamma(m,n)$ if and only if $\alpha\in S(m,\Z_+)$.
\endproof

We shall now define analogs of the elements $c_n^{\M}$.
First, for any subset $Q\subseteq \N_{mn}$ and any partition
$\M=(M_1,\ldots,M_r)$ such that $|\M|=|\ts Q|$ we set
\beq\label{cmq}
c_Q^{\M}=\sum (i_1,\ldots,i_{M_1})(j_1,\ldots,j_{M_2})\ldots
(k_1,\ldots,k_{M_r}),
\end{equation}
where $(i_1,\ldots,i_{M_1})$ etc. are cyclic permutations of the
corresponding indices and the summation is taken over all the orderings
$(i_1,\ldots,i_{M_1};j_1,\ldots,j_{M_2};\ldots;k_1,\ldots,k_{M_r})$ of the
elements of $Q$. We shall suppose that $c_{\emptyset}^{\emptyset}=1$.

Second, for any $\alpha\in\Gamma(m,\Z_+)$ and any subset $P\subseteq
\N_{mn}$ such that $\ord\alpha=|\ts P|$, we set
\beq\label{salp}
\Gamma(\alpha,P)=\{\sigma\in \Gamma(m,n)\ |\ \alpha(\sigma)=\alpha,
\ P(\sigma)=P,\ \M(\sigma)=(1^{n-m-|\ts P|}) \}, 
\end{equation}
i.e., $\sigma$ has to fix all the points in $\N_{mn}\setminus P$.

\bde\label{def:calpham}{\rm 
For any couple $(\alpha,\M)$, where
$\alpha\in\Gamma(m,\Z_+)$ and $\M$ is a partition such that
$\ord\alpha+|\M|\leq n-m$ we set
\beq\label{calpm}
c_n^{\alpha,\M}=\sum_{P,\ts Q}\sum_{\sigma\in \Gamma(\alpha,P)} 
\sigma\ts c_Q^{\M},
\end{equation}
where $P, Q$ are disjoint subsets in $\N_{mn}$ such that
\beq\label{pq}
|\ts P|=\ord\alpha,\qquad |\ts Q|=|\M|.
\end{equation}
}\endproof
\ede

\bpr\label{prop:zmnbas} Each of the families
\beq\label{zmnbas1}
c_n^{\alpha,\M}\qquad\text{with}\quad
\ord\alpha+|\M|=n-m,
\end{equation}
and
\beq\label{zmnbas2}
c_n^{\alpha,\M}\qquad\text{with}\quad
\ord\alpha+|\M|\leq n-m,\quad\text{and $\M$ has no part equal to $1$},
\end{equation}
forms a basis of $Z_m(n)$.
\epr

\Proof Note that $c_n^{\alpha,\M\cup 1\cup\ldots\cup 1}$ is proportional to
$c_n^{\alpha,\M}$. Therefore, it suffices to consider the family
\eqref{zmnbas1}.
By Proposition~\ref{prop:param}
the elements $c_n^{\alpha,\M}$ with
$\ord\alpha+|\M|=n-m$ are proportional to characteristic functions
of the $S_m(n)$-orbits in $\Gamma(m,n)$.
\endproof

Now we introduce analogs of the elements $\Delta_n^{\M}$.

\bde\label{def:delanm}{\rm 
For any couple $(\alpha,\M)$, where
$\alpha\in\Gamma(m,\Z_+)$ and $\M$ is a partition such that
$\ord\alpha+|\M|\leq n-m$ we set
\beq\label{delanm}
\Delta_n^{\alpha,\M}=
\sum_{P,\ts Q}\sum_{\sigma\in \Gamma(\alpha,P)}\epsi(P)\ts\sigma\ts
c_Q^{\M}\ts\epsi(Q)\ts\epsi(P), 
\end{equation}
where $\epsi(I):=(1-\ve_{i_1})
\cdots (1-\ve_{i_k})$ for $I=\{i_1,\dots,i_k\}$.
Here $P, Q$ are disjoint subsets in $\N_{mn}$ 
satisfying \eqref{pq}. We set $\Delta_n^{1,\emptyset}=1$,
where $\emptyset$ stands for the empty partition.
}\endproof
\ede

Note that \eqref{delanm} can be written in an equivalent form 
where the term $\epsi(Q)$ takes the leftmost position;
cf. \eqref{deltanm} and \eqref{del1}.

\bpr\label{prop:delinanm} The elements $\Delta_n^{\alpha,\M}$ 
belong to the algebra $A_m(n)$.
\epr

\Proof The semigroup
$\Gamma_m(n)$ is generated by the subgroup $S_m(n)$ and the idempotents
$\epsi_{m+1},\ldots, \epsi_n$. Therefore, it suffices to check that
$\Delta_n^{\alpha,\M}$ is stable under the conjugation by the
elements of $S_m(n)$ and commutes with the idempotents. 
The first claim is immediate from \eqref{delanm}. 
The second claim is verified exactly as its counterpart
for the elements 
$\Delta_n^{\M}$; see the proof of Proposition~\ref{prop:delao}.
\endproof

The following is an analog of Proposition~\ref{prop:thetdel}
and it is proved by the same argument.

\bpr\label{prop:thdelam} 
We have
\beq\label{thdelam}
\theta_{n}(\Delta_n^{\alpha,\M})=\Delta_{n-1}^{\alpha,\M}, 
\end{equation}
where we adopt the convention that
\beq\label{convam}
\Delta_k^{\alpha,\M}=0\quad \text{if\quad $ \ord\alpha+|\M|>k-m$.}
\end{equation}
\endproof
\epr

Our aim now is to prove an analog of Propositions~\ref{prop:bascent}
and \ref{prop:basao}; see Proposition~\ref{prop:basamn} below.
We need the following three lemmas.

\ble\label{lem:inamn} For $m<n$
\beq\label{inamn}
I(n)\cap
A_m(n)=(1-\epsi_{m+1})\cdots(1-\epsi_n)\ts
Z_m(n)\ts (1-\epsi_{m+1})\cdots(1-\epsi_n).
\end{equation}
\ele

\Proof Suppose that $x\in A(n)$ can be written as
\beq\label{xepy}
x=(1-\epsi_{m+1})\cdots(1-\epsi_n)\ts y\ts 
(1-\epsi_{m+1})\cdots(1-\epsi_n),
\end{equation}
where $y\in Z_m(n)$. Then $x\in A_m(n)$ since $x$ is invariant under
the conjugation by the elements of $S_m(n)$ 
and is annihilated when multiplied (from
the left or from the right) by any idempotent $\epsi_{m+1},\ldots,\epsi_n$.
Moreover, this also implies that $x\in I(n)$. 

Conversely, suppose $x\in I(n)\cap A_m(n)$. Then 
$x\ts \epsi_n=\epsi_n\ts x=0$.
Using the invariance of $x$ under the conjugation by the elements of
$S_m(n)$ we obtain $x\ts\epsi_i=\epsi_i\ts x=0$ for $i=m+1,\ldots,n$. 
Thus $x$ is invariant under the multiplication by
$(1-\epsi_{m+1})\cdots(1-\epsi_n)$ both from the left and from the right.

Further, we can write $x=y+y'$ where $y$ and $y'$ are spanned by
elements of $\Gamma(m,n)$ and $\Gamma(n)\setminus \Gamma(m,n)$, respectively. 
However,
\beq\label{eyprim}
(1-\epsi_{m+1})\cdots(1-\epsi_n)\ts y'\ts(1-\epsi_{m+1})\cdots(1-\epsi_n)=0
\end{equation}
since for each element $\gamma\in \Gamma(n)\setminus \Gamma(m,n)$
there exists $i>m$ such that $\gamma\ts\epsi_i=\gamma$ or 
$\epsi_i\ts\gamma=\gamma$.
This implies
\beq\label{xeye}
x=(1-\epsi_{m+1})\cdots(1-\epsi_n)\ts y\ts(1-\epsi_{m+1})\cdots(1-\epsi_n).
\end{equation}
Finally, averaging over the group $S_m(n)$ transforms $y$
into an element of $Z_m(n)$; cf. the proof of Lemma~\ref{lem:inao}.
\endproof

For $\sigma\in\Gamma(n)$, set
\beq\label{qsig}
Q(\sigma)=\{i\in\N_{mn}\ |\ \sigma_{ii}=1\}.
\end{equation}

\ble\label{lem:bij} 
The mapping
\beq\label{sigtogam}
\sigma\mapsto\gamma,\qquad
\gamma=\sigma\ts\epsi^{}_{Q(\sigma)}=\epsi^{}_{Q(\sigma)}\ts\sigma 
\end{equation}
defines a bijection of $\Gamma(m,n)$ onto the set of all $\gamma\in\Gamma(n)$
satisfying the conditions
\begin{align}\label{dor}
\dom\gamma\cap \N_{mn}&=\range\gamma\cap \N_{mn},\\
\label{degmg}
\deg_m\gamma&=n-m. 
\end{align}
\ele

\Proof The effect of the multiplication of $\sigma$ by $\epsi^{}_{Q(\sigma)}$
from the left or from the right consists of replacing all the diagonal
entries $\sigma_{ii}=1$ with $i>m$ by zeros. 
Therefore $\gamma$ satisfies \eqref{dor}.
Relation \eqref{degmg} follows from this observation and the fact that
both $\dom\sigma$ and $\range\sigma$ contain $\N_{mn}$.

Conversely, let $\gamma\in\Gamma(n)$ satisfy \eqref{dor} and \eqref{degmg}. 
Note that
\eqref{dor} can be reformulated as follows: for any $i=m+1,\ldots,n$ 
the $i$-th
row and the $i$-th column are zero or non zero at the same time, whereas
\eqref{degmg} means that all the diagonal entries 
$\gamma_{ii}$ with $i>m+1$ vanish.
Now, let $\sigma$ be defined by
\begin{align}
\sigma_{ij}&=\gamma_{ij}\qquad &{}&\text{if either $i\ne j$ or
${\rm{min}}\{i,j\}\leq m$},
\non\\
\sigma_{ii}&=1\qquad &{}&\text{if $i\in\ts \N_{mn}$}
\non
\end{align}
and the $i$-th row (or
the $i$-th column) of $\gamma$ is zero.
Then it is easy to see that $\sigma\in \Gamma(m,n)$ and that $\gamma$ is the
image of $\sigma$ under the mapping \eqref{sigtogam}. 
\endproof

\ble\label{lem:projamn} For $m<n$ the restriction of the projection 
$\theta_{n}:A_m(n)\to A_m(n-1)$ to the subspace
$F_m^{n-m-1}(A_m(n))$ is injective.
\ele

\Proof Suppose that $x\in A_m(n)$ and $\theta_{n}(x)=0$. We will show
that then $\deg_m x=n-m$ unless $x=0$.

By Lemma~\ref{lem:inamn}, $x$ can be written as a linear combination of the
elements
of type
\begin{multline}\label{epsig}
(1-\epsi_{m+1})\cdots(1-\epsi_n)\ts\sigma
\ts(1-\epsi_{m+1})\cdots(1-\epsi_n)\\
=\sum_{R,S\subseteq N_{mn}}(-1)^{|\ts R|+|\ts S|} \epsi^{}_R\ts\sigma
\ts\epsi^{}_S, \qquad\sigma\in \Gamma(m,n).	           
\end{multline}
Let us divide the terms in the sum \eqref{epsig} into two groups depending on
whether $R\cup S$ 
contains $Q(\sigma)$ or not. Then the terms
of the first group are of $m$-degree 
$n-m$ whereas those of the
second group are of $m$-degree ${}<n-m$. So, it suffices to prove that
the elements
\beq\label{gro}
\sum_{R\cup S\supseteq\ts
Q(\sigma)}(-1)^{|\ts R|+|\ts S|}\epsi^{}_R\ts\sigma
\ts\epsi^{}_S, \qquad \sigma\in \Gamma(m,n),
\end{equation}
are linearly independent.
Note that, in the case $R\cup S=Q(\sigma)$,
\beq\label{casex1}
\epsi^{}_R\ts\sigma
\ts\epsi^{}_S=\sigma\ts\epsi^{}_{Q(\sigma)}\qquad\text{and}\qquad
\rank\sigma\ts\epsi^{}_{Q(\sigma)}=\rank\sigma-|\ts Q(\sigma)|,
\end{equation}
whereas, in the case $R\cup S$ strictly contains $Q(\sigma)$,
\beq\label{casex2}
\rank\sigma\ts\epsi^{}_{Q(\sigma)}<\rank\sigma-|\ts Q(\sigma)|.
\end{equation}
Therefore, we now need to show
that for any fixed $k$ the elements
\beq\label{linindq}
\sum_{R\cup
S=Q(\sigma)}(-1)^{|\ts R|+|\ts S|}\sigma\ts\epsi^{}_{Q(\sigma)}, 
\end{equation} 
where $\sigma$ runs over the subset of the elements in $\Gamma(m,n)$ with
$\rank\sigma-|\ts Q(\sigma)|=k$, are linearly independent.

Lemma~\ref{lem:bij} implies that the elements
$\sigma\ts\epsi^{}_{Q(\sigma)}\in\Gamma(n)$ are pairwise distinct.
Hence it remains to
prove that all the coefficients in \eqref{linindq} 
are non-vanishing. This is implied by the 
following general fact: if $Q$ is an arbitrary
finite set, then
\beq\label{genfa}
\sum_{R,S\subseteq Q,\ R\cup S=Q}(-1)^{|\ts R|+|\ts S|}\ne 0.
\end{equation}
We will prove that the sum in \eqref{genfa} equals $(-1)^q$ where $q=|\ts Q|$.
Indeed, for any $r=0,1,\ldots,q$, there are 
$\displaystyle{\binom nr}$ subsets
$R\subseteq Q$ with $|\ts R|=r$. Given $R$, for any
$t=0,1,\ldots,r$, there are 
$\displaystyle{\binom rt}$ subsets $S\subseteq Q$ such that
$R\cup S=Q$ and $|\ts R\cap S|=t$. Since
\beq\label{cardd}
|\ts R|+|\ts S|=r+t+(q-r)=q+t,
\end{equation}
the sum in \eqref{genfa} equals
\beq
(-1)^q\sum_{r=0}^q \binom nr \sum_{t=0}^r(-1)^t\binom rt . 
\non
\end{equation}
If $r=0$ then the interior sum is equal to $1$,
otherwise it is zero. Therefore the entire sum is $(-1)^q$.
\endproof

\bpr\label{prop:basamn} The elements $\Delta_n^{\alpha,\M}$ with
\beq\label{ordin}
\ord\alpha+|\M|\leq n-m
\end{equation}
form a basis of $A_m(n)$. Moreover, for any $M$ with 
$0\leq M\leq n-m$
the elements $\Delta_n^{\alpha,\M}$ satisfying
\beq\label{ordine}
\ord\alpha+|\M|\leq M
\end{equation}
form a basis of $F_m^M(A_m(n))$. 
\epr

\Proof It suffices to prove the second claim.
We use induction on $n$ and follow the argument of the
proof of Proposition~\ref{prop:basao}.
The claim is obviously true for $n=m$. 
Assume that $n\geq m+1$ 
and $M\leq n-m-1$. Lemma~\ref{lem:projamn} implies that 
the elements $\Delta_n^{\alpha,\M}$
with $\ord\alpha+\vert\M\vert\leq M$ form a basis of $F_m^M(A_m(n))$.

To show that the elements $\Delta_n^{\alpha,\M}$ 
with $\ord\alpha+\vert\M\vert=n-m$ 
form a basis of $I(n)\cap A_m(n)$ note that
\beq\label{delacam}
\Delta_n^{\alpha,\M}=(1-\epsi_{m+1})\cdots (1-\epsi_n)\ts
c_n^{\alpha,\M}\ts(1-\epsi_{m+1})\cdots (1-\epsi_n); 
\end{equation}
see \eqref{delanm}.
Now the claim follows from
Proposition~\ref{prop:zmnbas} and the fact 
that the elements $c_n^{\alpha,\M}$, being multiplied by
$(1-\epsi_{m+1})\cdots (1-\epsi_n)$, remain linearly independent;
cf. \eqref{compos}. 
\endproof

Using Proposition~\ref{prop:thdelam} we can introduce the elements
$\Delta^{\alpha,\M}\in A_m$ as sequences
$\Delta^{\alpha,\M}=(\Delta_n^{\alpha,\M}\ |\ n\geq m)$.

\medskip
\noindent
{\it Remark\/.} We can regard  $\Delta^{\alpha,\M}$ as a formal series
given by \eqref{delanm} where the sum is taken over all disjoint subsets
$P$ and $Q$ in $\{m+1,m+2,\dots\}$ satisfying \eqref{pq}. \endproof

\bth\label{thm:basam} The elements $\Delta^{\alpha,\M}$ with
$\alpha\in \Gamma(m,\Z_+)$ and $\M\in\PP$
form a basis of the algebra $A_m$. Moreover, for any $M\geq 0$,
the elements $\Delta^{\alpha,\M}$ 
with $\ord\alpha+\vert\M\vert\leq M$ form a basis of the
$M$-th subspace $F_m^M(A_m)$ in $A_m$.
\eth

\Proof The first claim follows from the second one. The second
claim follows from Proposition~\ref{prop:basamn}
and the definition of $F_m^M(A_m)$ as the projective limit of the
spaces $F_m^M(A_m(n))$.
\endproof

\bco\label{cor:isomam} For $n>M$, the mapping
\beq\label{isomamn}
\theta_{n}: F_m^M(A_m(n))\to F_m^M(A_m(n-1))
\end{equation}
is an isomorphism of vector spaces and so is
the mapping
\beq\label{isomam}
\theta^{(n)}:F_m^M(A_m)\to F_m^M(A_m(n)), \qquad	 n\geq M.
\end{equation}
In particular, $\dim F_m^M(A_m)<\infty$.	
\endproof
\eco

For each $k=1,\dots,m$ consider
the following elements of $A_m(n)$
\beq\label{jmre}
u_{k|n}=\sum_{i=k+1}^n(ki)(1-\ve_k)(1-\ve_i)=
\sum_{i=k+1}^n(1-\ve_i)(ki)(1-\ve_i).
\end{equation}
The image of $u_{k|n}$ under the retraction homomorphism \eqref{retr}
is the {\it Jucys--Murphy element\/} for $S(n)$; cf. \cite{j:sp},
\cite{m:nc}.
We obviously have $\theta_n(u_{k|n})=u_{k|n-1}$ and so,
for each $k$ the element $u_k\in A_m$ can be defined as
the sequence $u_k=(	u_{k|n}\ |\ n\geq m)$. Recall that the algebra
$A(m)$ is naturally embedded in $A_m$; see Proposition~\ref{prop:stabemb}.

\bpr\label{prop:ukcom}
The following relations hold in the algebra $A_m$:
\begin{align}\label{defasn}
s_k\ts u_k&=u_{k+1}\ts s_k+(1-\ve_k)(1-\ve_{k+1}),
\qquad &s_k\ts u_l&=u_l\ts s_k,\quad l\ne k,k+1;\\
\label{defasn1}
u_k\ts u_l&=u_l\ts u_k,\qquad\qquad \ve_k\ts u_k=u_k\ts \ve_k=0,\qquad
&\ve_i\ts u_k&=u_k\ts \ve_i,\quad i\ne k;
\end{align}
where $s_k=(k,k+1)$.
\epr

\Proof For $n>m$ we have $u_{1|n}=\Delta^{(2)}_n-\Delta^{\prime(2)}_{n-1}$,
where $\Delta^{\prime(2)}_{n-1}$ is the element
of the center of $\C[\Gamma_1(n)]$ given by \eqref{deltanm} with
the sum taken over the indices from $\{2,\dots,n\}$.
Now an easy induction proves that the elements $u_{1|n},\dots,u_{m|n}$
pairwise commute, and so do the elements $u_1,\dots,u_m$.
The remaining relations easily follow from \eqref{jmre} and
the relations in the
algebra $A(n)$. \endproof

We shall denote by  $\wH_m$ the subalgebra of 
$A_m$ generated by $A(m)$ and the elements $u_1,\dots,u_m$.
The following is our main result. The theorem describes
the structure of the algebra $A_m$.

\bth\label{thm:am}
We have an algebra isomorphism
\beq\label{am}
A_m\simeq A_0\ot \wH_m.	
\end{equation}
Moreover, the algebra $\wH_m$ is isomorphic to an abstract algebra
with generators \newline
$s_1,\dots,s_{m-1}$, $\ve_1,\dots,\ve_m$,
$u_1,\dots,u_m$ and the defining relations	
given by \eqref{defsn}--\eqref{defsga} and 
\eqref{defasn}--\eqref{defasn1}.
\eth

Recall that by Theorem~\ref{thm:alaisom} $A_0$ is 
isomorphic to the algebra of shifted symmetric functions $\Lambda^*$.

\Proof
For any $\alpha\in\Gamma(m,\Z_+)$ and $\M\in\PP$ such that
$\ord \alpha+|\M|\leq n-m$ we have the equality in the
algebra $A_m(n)$,
\beq\label{dadm}
\Delta^{\alpha,\emptyset}_n\Delta^{1,\M}_n=\Delta^{\alpha,\M}_n
+\text{\ lower $m$-degree terms},
\end{equation}
where $\emptyset$ stands for the empty partition while
$1\in \Gamma(m,\Z_+)$ is the $m\times m$ identity matrix; 
cf. the proofs of Proposition~\ref{prop:prodcn} 
and Corollary~\ref{cor:genaon}.
On the other hand, we have
\beq\label{delmdel}
\deg_m(\Delta^{\M}_n-\Delta^{1,\M}_n)< |\M|;
\end{equation}
see \eqref{deltanm}	for the definition of $\Delta^{\M}_n$.
Now \eqref{dadm} and
Proposition~\ref{prop:basamn} imply that the elements
$\Delta^{\alpha,\emptyset}_n\Delta^{\M}_n$ with
$\ord \alpha+|\M|\leq n-m$ form a basis	of $A_m(n)$.
Hence the elements $\Delta^{\alpha,\emptyset}\Delta^{\M}$
with $\alpha\in\Gamma(m,\Z_+)$ and $\M\in\PP$ form
a basis of the algebra $A_m$. In other words,
$A_m$ is a free $A_0$-module
with the basis $\{\Delta^{\alpha,\emptyset}\ |\ \alpha\in\Gamma(m,\Z_+)\}$.

Further, if 
$\alpha\in \Gamma(m,\Z_+)$ has zero rows
$i_1,\dots,i_r$ then 
\beq\label{alepr}
\Delta_n^{\alpha,\emptyset}=\ve_{i_1}\cdots 
\ve_{i_r}\ts\Delta_n^{\alpha',\emptyset}
\end{equation}
for some element $\alpha'\in S(m,\Z_+)$. Observe now
that every element
$\alpha'\in S(m,\Z_+)$ 	can be written as a product of the form
\beq\label{alsig}
\alpha'=\sigma\ts \alpha_1^{k_1}\cdots \alpha_m^{k_m},\qquad k_i\geq 0,
\end{equation}
where $\sigma\in S(m)$ and $\alpha_i\in\Gamma(m,\Z_+)$ 
is the diagonal matrix whose
$ii$-th entry is $z$ and all other diagonal entries are equal to $1$.
This implies that modulo lower $m$-degree
terms, the element $\Delta_n^{\alpha',\emptyset}$ coincides
with the product 
\beq\label{prosig}
\Delta_n^{\sigma,\emptyset}(\Delta_n^{\alpha_1,\emptyset})^{k_1}
\cdots (\Delta_n^{\alpha_m,\emptyset})^{k_m};
\end{equation}
cf. the proof of \eqref{unipart}. The claim remains valid if
we replace each $\Delta_n^{\alpha_k,\emptyset}$ with the element
$u_{k|n}$. Indeed, this follows from the equality
\beq\label{ukdel}
\Delta^{\alpha_k,\emptyset}_n=u_{k|n} +\text{\ elements of $m$-degree zero}.
\end{equation}
Note also that the element
$\Delta_n^{\sigma,\emptyset}$ can be identified with $\sigma$.
Thus, modulo lower $m$-degree terms, 
the element \eqref{prosig}	coincides with
the product $\sigma\ts u^{k_1}_{1|n}\cdots u^{k_m}_{m|n}$.
Using an obvious induction on the $m$-degree
we may conclude that the $A_0$-module $A_m$ is 
generated by the subspace $\wH_m$.

To prove that $\wH_m$ generates the $A_0$-module $A_m$ freely,
we check that for any $M>0$ the dimension of the subspace 
$F^M_m(\wH_m)$ is less or equal to the number of elements
$\alpha\in\Gamma(m,\Z_+)$ with $\ord\alpha\leq M$.
Indeed, by Proposition~\ref{prop:ukcom} the subalgebra
$\wH_m$ is spanned by the elements of the form
$\gamma\ts u_1^{k_1}\cdots u_m^{k_m}$ with $\gamma\in\Gamma(m)$.
The relation $\ve_k\ts u_k =0$ ensures that
such a product  
is zero unless $k_j=0$ for
each zero column $j$ in $\gamma$. 
To each of the 
nonzero products associate the element $\alpha\in \Gamma(m,\Z_+)$
which has the $ij$-entry $z^{k_j}$ where
the $j$-th column of $\gamma$ is nonzero with $\gamma_{ij}=1$.
This shows that the cardinality
of a basis of $F^M_m(\wH_m)$ can be
at most the number of elements $\alpha\in \Gamma(m,\Z_+)$
with $\ord\alpha\leq M$, proving \eqref{am}.

To prove the second claim of the theorem note that by
Proposition~\ref{prop:ukcom} there is an algebra epimorphism
from the abstract algebra in question to $\wH_m$. The above argument
implies that the nonzero products 
$\gamma\ts u_1^{k_1}\cdots u_m^{k_m}$ with $\gamma\in\Gamma(m)$
form a basis of $\wH_m$.
\endproof

\bco\label{cor:amngen}
The mapping
\beq\label{amngen}
s_k\mapsto s_k,\qquad \ve_k	\mapsto \ve_k,\qquad u_k\mapsto u_{k|n}
\end{equation}
defines an algebra homomorphism $\psi:\wH_m\to A_m(n)$. The algebra
$A_m(n)$ is generated by $A_0(n)$ and the image of $\psi$.
\endproof
\eco

The {\it degenerate affine Hecke algebra\/}	$\HH_m$	(see \cite{d:da},
\cite{l:ah})
is defined to be generated by elements $s_1,\dots,s_{m-1}$ and
$u_1,\dots,u_m$ with the defining relations \eqref{defsn}
and 
\begin{align}\label{defahe}
s_k\ts u_k&=u_{k+1}\ts s_k+1,
\qquad &s_k\ts u_l&=u_l\ts s_k,\quad l\ne k,k+1;\\
\label{defahe1}
u_k\ts u_l&=u_l\ts u_k. &&
\end{align}
As a linear space, $\HH_m$ is isomorphic to the tensor product
$\C[S(m)]\ot\C[u_1,\dots,u_m]$.	The following corollary
is implied by Theorem~\ref{thm:am} and provides an analog of
the retraction homomorphism \eqref{retr}.

\bco\label{cor:homdegen}
The mapping
\beq\label{homdegen}
s_k\mapsto s_k,\qquad u_k\mapsto u_k,\qquad \ve_k\to 0
\end{equation}
defines an algebra epimorphism $\wH_m\to \HH_m$.
\endproof
\eco

It can be seen from
the proof of Theorem~\ref{thm:am}
that the retraction homomorphisms \eqref{retr} and \eqref{homdegen}
``respect" the homomorphism $\psi:\wH_m\to A_m(n)$
defined in Corollary~\ref{cor:amngen}. 
More precisely, the following
result takes place.
It was announced in \cite[Theorem~11]{o:ea},
and a proof was given in \cite{ov:na}. We denote
by $B_m(n)$ the centralizer of $S_m(n)$ in the group
algebra $\C[S(n)]$; see Introduction.

\bco\label{cor:amnhec}
The mapping
\beq\label{amnhec}
s_k\mapsto s_k,\qquad u_k\mapsto \sum_{i=k+1}^{n}(ki)
\end{equation}
defines an algebra homomorphism $\varphi:\HH_m\to B_m(n)$. The algebra
$B_m(n)$ is generated by $B_0(n)$ and the image of $\varphi$.
\endproof
\eco

\end{document}